\documentclass{amsart}

\usepackage{amssymb,amsmath,amsthm,texdraw}
\usepackage{hyperref}
\usepackage{color}

\newtheorem{Theorem}{Theorem}
\newtheorem{Proposition}{Proposition}
\newtheorem{Lemma}{Lemma}
\newtheorem{Definition}{Definition}
\newtheorem{Corollary}{Corollary}

\newtheorem{Theo}{Theorem}

\newtheorem{Prop}{Proposition}

\theoremstyle{remark}
\newtheorem{Remark}{Remark}
\newtheorem{Example}{Example}
\newtheorem{Notation}{Notation}

\title{Infinitesimal adjunction and polar curves}

\author{Nuria Corral}

\address{Dpto. Matem\'{a}ticas, Estad\'{\i}stica y Computaci\'{o}n. Universidad
de Cantabria. Avda. de los Castros s/n, 39005 Santander. Spain}
\email{nuria.corral@unican.es}

\thanks{The author was partially supported by the research projects
MTM2007-66262 (Ministerio de Educaci\'{o}n y Ciencia),
MTM2006-15338-C02-02 (Ministerio de Educaci\'{o}n y Ciencia), VA059A07
(Junta de Castilla y Le\'{o}n) and PGIDITI06PXIB377128PR (Xunta de
Galicia).}



\begin{document}
\begin{abstract}
The polar curves of foliations $\mathcal F$ having a curve $C$ of
separatrices generalize the classical polar curves associated to
hamiltonian foliations of $C$. As in the classical theory, the
equisingularity type ${\wp}({\mathcal F})$ of a generic polar curve
depends on the analytical type of ${\mathcal F}$, and hence of $C$.
In this paper we find the equisingularity types $\epsilon (C)$ of
$C$, that we call kind singularities, such that ${\wp}({\mathcal
F})$ is completely determined by $\epsilon (C)$ for Zariski-general
foliations $\mathcal F$. Our proofs are mainly based on the
adjunction properties of the polar curves. The foliation-like
framework is necessary,  otherwise we do not get the right concept
of general foliation in Zariski sense and, as we show by examples,
the hamiltonian case can be out of the set of general foliations.
\end{abstract}

\maketitle

\section{Introduction}
Let $\mathcal F$ be a germ of holomorphic foliation of $ ({\mathbb
C}^2,0)$ having a curve of sepa\-ratrices $C$. The {\em polar
curve\/} $\Gamma$ of $\mathcal F$ with respect to a direction
$[a:b]\in {\mathbb P}^1_{\mathbb C}$ is given by  $\omega \wedge (a
dy-bdx)=0$, where $\omega$ is a $1$-form defining $\mathcal F$.
There is a Zariski-open set of directions such that the
equisingularity type  $\epsilon (\Gamma\cup C)$ of $\Gamma\cup C$ is
the same one, independent of $\omega$ and of the coordinates. We
denote $\wp({\mathcal F})$ this generic type of equisingularity.
This paper is devoted to provide an accurate description of the
types $\wp({\mathcal F})$ in terms of the equisingularity type
$\epsilon(C)$ of $C$.

We work with foliations in the class ${\mathbb G}_{C}^{\ast}$ of the
generalized curves without ``bad resonances'' defined as follows. A
foliation ${\mathcal F}$ belongs to ${\mathbb G}_C^{\ast}$ if
\begin{enumerate}
\item It is a generalized curve in the sense of Camacho-Lins Neto-Sad
(\cite{Cam-S-LN}) having $C$ as curve of separatrices. Note that, in
this case, the minimal morphism of reduction of singularities
$\pi_C$ of $C$ is also the reduction of singularities of ${\mathcal
F}$.
\item For any $C$-{\em ramification\/} $\rho:  ({\mathbb C}^2,0)
\rightarrow ({\mathbb C}^2,0)$ (that is, $\rho$ is transversal to
$C$ and $\rho^{-1}C$ has only non-singular branches), there is no
corner in the reduction of singularities of $\rho^{\ast}{\mathcal
F}$ with Camacho-Sad index equal to $-1$.
\end{enumerate}
If $C=(f=0)$, the hamiltonian foliation $df=0$ belongs to  ${\mathbb
G}_C^{\ast}$. But the class ${\mathbb G}_C^{\ast}$ is wider than
that.  Let us write $f=\prod_{i=1}^r f_i$, then the logarithmic
foliations
$${\mathcal L}_\lambda=\left(
\sum_{i=1}^r\lambda_i\frac{df_i}{f_i}=0\right)
$$
belong to this class if $\lambda=(\lambda_1,\cdots,\lambda_r)$ avoid
certain rational resonances. More generally, each generalized curve
foliation $\mathcal F$  has a well defined {\em logarithmic model\/}
${\mathcal L}_{\lambda}$, $\lambda=\lambda({\mathcal F})$, of the
above type such that the Camacho-Sad indices of $\mathcal F$ and
${\mathcal L}_{\lambda}$ coincides along the reduction of
singularities \cite{Cor-03}.

There is a first relationship between $\epsilon(C)$ and
$\wp({\mathcal F})$ described in the {\em decomposition theorem of
the polar curve\/} \cite{Cor-03}, proved  by several authors in
different contexts \cite{Mer,Kuo-L,Gar,Rou}. It can be stated as
follows:

\begin{Theo}[Decomposition \cite{Cor-03}] Let $\rho$ be a $C$-ramification.
If $\Gamma$ is a generic polar curve of ${\mathcal F} \in {\mathbb
G}_{C}^{\ast}$, then $\rho^{-1}\Gamma$ is a strict adjoint of
$\rho^{-1}C$.
\end{Theo}
If  $Y\subset ({\mathbb C}^2,0)$ is a curve with only non-singular
branches, we say that a curve $Z \subset ({\mathbb C}^2,0)$ is a
{\em strict adjoint\/} of $Y$ if  the multiplicities satisfy
$m_p(Z)=m_p(Y)-1$ at the infinitely near points $p$ of $Y$ and $Z$
does not go through the corners of the desingularization of $Y$.
(Compare with the definition in \cite{Cas-00}, p. 152).

There are infinitely many possible equisingularity types $\epsilon
(Y\cup Z)$ for a fixed $Y$ and $Z$ being strict adjoint of $Y$. In
section~\ref{sec:loc-inv}  we prove the following result of
finiteness by using a control of the Newton polygon of a generic
polar curve $\Gamma$ (a similar result for the case of hamiltonian
foliations can be deduced from the virtual behaviour of the polar
curves described in \cite{Cas-00}).
\begin{Theo}
There exists a finite number of equisingularity types $\wp({\mathcal
F})$, where ${\mathcal F} \in {\mathbb G}_{C'}^{\ast}$ and $C'$ is
such that $\epsilon(C')=\epsilon(C)$.
\end{Theo}

Take as above $Y\subset ({\mathbb C}^2,0)$  with only non-singular
branches. A strict adjoint curve $Z$ of $Y$ is a {\em perfect
adjoint\/} curve of $Y$ if $\pi_Y$ desingularizes $Z$. In this case
the equisingularity type $\chi_Y=\epsilon (Y \cup Z)$ does not
depend on $Z$. Section~\ref{sec:non-sp} is devoted to prove the
following result of genericity
\begin{Theo}[of genericity]\label{th:genericity} Assume that $C$ has only
non-singular branches. There is a non-empty Zariski-open set $U_C
\subset {\mathbb P}_{\mathbb C}^{r-1}$ defined by

\begin{quote}
``$\lambda \in U_C$ if there exists ${\mathcal F} \in {\mathbb
G}_{C}^{\ast}$ with $\wp({\mathcal F})=\chi_C$ and
$\lambda=\lambda({\mathcal F})$".
\end{quote}
 Moreover, for each ${\mathcal F}
\in {\mathbb G}_{C}^{\ast}$ with $\lambda({\mathcal F}) \in U_C$ we
have that $\wp ({\mathcal F})=\chi_C$.
\end{Theo}
In general, it is not possible to define $\chi_C$ in a  way
compatible with $C$-ramifications. This is the characteristic
property of the {\em kind equisingularity types\/} that we introduce
below.

Let $G(C)$ be the dual graph of $C$ oriented by its first divisor.
Associate to each divisor $E$ the multiplicity $m(E)$  given by any
$E$-``curvette" and the number $b_E$ of  edges and arrows which
leave from $E$. Thus  $E$ is a {\em bifurcation divisor\/} if $b_E
\geq 2$ and a {\em terminal divisor\/} if $b_E=0$. A {\em dead
arc\/} joins a bifurcation divisor with a terminal divisor, with no
other bifurcations. We say that $\epsilon(C)$ is {\em kind \/} if
$m(E_b)=2m(E_t)$, for each dead arc of $G(C)$ starting at $E_b$ and
ending at $E_t$. The next proposition, proved in
section~\ref{sec:tipo-amable}, gives a characterization of kind
equisingularity types in terms of adjunction
\begin{Prop} The equisingularity type
$\epsilon (C)$ is kind if and only if there is a germ of curve $Z
\subset ({\mathbb C}^2,0)$ such that $\rho^{-1}Z$ is a perfect
adjoint of $\rho^{-1}C$ for any $C$-ramification $\rho$. Moreover
$\epsilon(C\cup Z)$ does not depend on the choice of $Z$.
\end{Prop}
For kind equisingularity types we define $\chi_C=\epsilon(C\cup Z)$
and we say that such $Z$ are perfect adjoint curves of $C$. The next
proposition, proved in section~\ref{sec:tipo-amable}, gives a
precise description of $\chi_C$ for kind equisingularity types. (For
classical polar curves, our description is slightly more precise
than the one in \cite{Le-M-W}).
\begin{Prop}
Let $C$ be a curve with kind equisingularity type and $Z$ a perfect
adjoint curve of $C$. Then $\pi_C$ gives a reduction of
singularities of $Z\cup C$. Moreover, the branches of $Z$ intersect
an irreducible component $E$ of the exceptional divisor of $\pi_C$
as follows:
\begin{itemize}
    \item If $E$ is a bifurcation divisor of $G(C)$,  the number
    of branches of $Z$ cutting $E$
    equals to $b_E-2$ if $E$ is in a dead arc and to
    $b_E-1$ otherwise.
    \item If $E$ is a terminal divisor of a dead arc of $G(C)$, there is
    exactly one branch of $Z$
    through $E$.
\item Otherwise, no branches of $Z$  intersect $E$.
\end{itemize}
\end{Prop}
Finally, in section~\ref{sec:main-th}, we relate the polar curves to
the adjoint curves in the case of kind equisingularity types. As a
consequence we obtain a precise description of  $\wp({\mathcal F})$
if $\epsilon (C)$ is kind. Let us define the Zariski open set $U_C
\subset {\mathbb P}_{\mathbb C}^{r-1}$ by
\begin{quote}
``$\lambda \in U_C$ if there exists ${\mathcal F} \in {\mathbb
G}_{C}^{\ast}$ with $\lambda=\lambda({\mathcal F})$ having a generic
polar curve $\Gamma$ such that $\rho^{-1}\Gamma$ is a perfect
adjoint of $\rho^{-1}C$, for  any $C$-ramification $\rho$"
\end{quote}
Then we prove the following theorem
\begin{Theo}
The curve $C$ has a kind equisingularity type if and only if $U_C
\ne \emptyset$.  In this case $\wp({\mathcal F})=\chi_C$ for any
${\mathcal F} \in {\mathbb G}_{C}^{\ast}$ such that
$\lambda({\mathcal F})\in U_C$.
\end{Theo}
The hamiltonian foliations $df=0$ have vector of exponents
$\lambda=\underline{1}$. We provide examples such that
$\underline{1}\notin U_C$, hence the consideration of the class
${\mathbb G}_{C}^{\ast}$ is essential for this theory.

The main results of this paper were announced in \cite{Cor-06}. Our
results are of local nature in the framework of foliations (see also
\cite{Rou,Cor-03,Cor-08}). The classical local study of polar curves
has been developed by several authors
(\cite{Tei,Mer,Kuo-L,Le-M-W,Cas-00,Gar}). There are also related
works for foliations from the global viewpoint \cite{Mol,Fas-P}.

\section{Strict adjoint curves}\label{sec:adjuntos}
Before starting the study of polar curves, we describe some
properties that can be deduced from the fact that a curve is a
strict adjoint of another curve. We recall the notion of a strict
adjoint curve:
\begin{Definition}
Assume that $C$ has only non-singular branches. We say that $Z$ is a
{\em strict adjoint\/} of $C$ if $m_p(Z) = m_p(C)-1$ at each
infinitely near point $p$ of $C$ and $Z$ does not go through the
corners of the desingularization of $C$.
\end{Definition}
If $Z$ is a strict adjoint of $C$, the properties above allow to
give a decomposition of $Z$ into bunches of branches in terms of the
equisingularity data of $C$. Let us describe it using the dual graph
$G(C)$ of $C$ which is constructed from the minimal reduction of
singularities $\pi_C : M \rightarrow ({\mathbb C}^2,0)$ of $C$ (see
appendix~\ref{ap:grafodual} for all the notations concerning the
dual graph of a curve). Given a divisor $E$ of $\pi_{C}^{-1}(0)$, we
denote by $\pi_E : M_E \rightarrow ({\mathbb C}^2,0)$ the morphism
reduction of $\pi_C$ to $E$ (see appendix~\ref{ap:grafodual});
recall that $\pi_C = \pi_E \circ \pi_{E}'$. Let $B(C)$ be the set of
bifurcation divisors of $G(C)$. For any $E \in B(C)$, we define
$Z^E$ to be the union of the branches $\zeta$ of $Z$ such that
\begin{itemize}
    \item $\pi_{E}^{\ast} \zeta \cap \pi_{E}^{\ast} C = \emptyset$
    \item If $E' < E$, then $\pi_{E}^{\ast} \zeta \cap
    \pi_{E}'(E')=\emptyset$
\end{itemize}
where $\pi_{E}^{\ast}\zeta$ denotes the strict transform of $\zeta$
by $\pi_{E}$. Thus there is a unique decomposition $Z = \cup_{E \in
B(C)} Z^E$ satisfying that:
\begin{enumerate}
    \item[d1.] $m_0(Z^E)=b_E-1$.
    \item[d2.] $\pi_{E}^* Z^{E} \cap \pi_{E}^* C = \emptyset$.
    \item[d3.] If $E' < E$ then $\pi_{E}^*Z^E \cap \pi_{E}' (E') =
    \emptyset$.
    \item[d4.] If $E' > E$ then $\pi_{E'}^* Z^{E} \cap E_{red}' =
    \emptyset$.
\end{enumerate}
In particular, if $E$ is not a bifurcation divisor we have that
$\pi_{E}^*Z \cap E_{red} = \pi_{E}^{\ast}C \cap E_{red}$. Moreover,
the properties  above imply the following ones which are stated in
terms of the coincidences and of the data in $G(C)$. For each
irreducible component $\zeta$ of $Z^E$ we have that
\begin{itemize}
    \item[(D-i)] ${\mathcal C}(C_i,\zeta) = v(E)$ if $E$ belongs to
    the geodesic of $C_i$;
    \item[(D-ii)] ${\mathcal C}(C_j,\zeta)={\mathcal C}(C_j,C_i)$
    if $E$ belongs to the geodesic of $C_i$ but not to the one of
    $C_j$.
\end{itemize}
(see appendix \ref{ap:grafodual} for the definitions of $b_E$,
$v(E)$ and {\em geodesic\/} of a curve in $G(C)$).

Consider now any curve $C$ and let $\rho: ({\mathbb C}^2,0)
\rightarrow ({\mathbb C}^2,0)$ be any $C$-ramification (the reader
can refer to appendix~\ref{ap:ramificacion} for notations and
general results concerning ramifications). If $\tilde{Z}=\rho^{-1}Z$
is a strict adjoint of $\tilde{C}=\rho^{-1}C$, then there is also a
decomposition of $Z$ in terms of the equisingularity data of $C$:
for any bifurcation divisor $E$ of $G(C)$, we define $Z^E$  to be
such that
$$
\rho^{-1} Z^E = \bigcup_{i=1}^{\underline{n}_E}
\tilde{Z}^{\tilde{E}^{j}},
$$
where $\{\tilde{E}^{j}\}_{j=1}^{\underline{n}_E}$ are the divisors
of $G(\tilde{C})$ associated to $E$ in $G(\tilde{C})$ and
$\tilde{Z}= \cup_{\tilde{E} \in G(\tilde{C})} \tilde{Z}^{\tilde{E}}$
is the decomposition of $\tilde{Z}$ described above. Hence, we get a
decomposition $Z = \cup_{E \in B(C)}Z^E$ such that:
\begin{itemize}
    \item[D1.] $m_0(Z^E)=\left\{
                           \begin{array}{ll}
                             \underline{n}_E n_E (b_E -1), & \hbox{if $E$ does not belong to a dead arc;} \\
                             \underline{n}_E n_E (b_E-1) - \underline{n}_E, & \hbox{otherwise.}
                           \end{array}
                         \right.$
    \item[D2.] $\pi_{E}^{\ast}Z^E \cap \pi_{E}^{\ast} C = \emptyset.$
    \item[D3.] If $E' < E$, then $\pi_{E}^{\ast} Z^E \cap \pi_{E}'
(E') = \emptyset$.
    \item[D4.] If $\pi_{E}^{\ast} Z^E \cap \pi_{E}'(E') \neq
\emptyset$, then $\pi_{E}'(E') > E_{red}$.
    \item[D5.] If $E'>E$ and $E'$ does not belong to a dead arc
joined to $E$, then $E'_{red} \cap \pi_{E'}^{\ast} Z^E =\emptyset$.
\end{itemize}
Moreover, properties (D-i) and (D-ii) also hold now for a branch
$\zeta$ of $Z^E$.

It is clear that the properties above do not determine the
equisingularity type of the curve $Z$ even if $C$ has only
non-singular branches. Let us introduce a definition:
\begin{Definition}
Assume that $C$ has only non-singular branches and let $Z$ be a
strict adjoint of $C$. We say that  $Z$  is a {\em perfect adjoint
curve\/} of $C$ if $\pi_C$ gives a reduction of singularities of
$Z$.
\end{Definition}
Let us state a criterion to check if a curve $Z$ is a perfect
adjoint of $C$.
\begin{Proposition}\label{criterio-adjunto-perfecto}
Let $C$ be a curve with only non-singular branches.  A strict
adjoint curve $Z$ of $C$ is perfect adjoint curve of $C$ if and only
if the set
$$\pi_{E}^*Z\cap E_{red} \smallsetminus \pi_{E}^* C \cap
E_{red}
$$
has exactly $b_E-1$ points for each irreducible component $E$ of
$\pi^{-1}_C(0)$.
\end{Proposition}
\begin{proof}
Observe that the second part of the statement always holds when $E$
is not a bifurcation divisor ($b_E=1$) since $\pi^{\ast}_E Z \cap
E_{red} = \pi^*_E C  \cap E_{red}$ (see the properties of the
decompositions  above). Therefore we only need to prove the result
for bifurcation divisors. Recall that there is a decomposition $Z =
\cup_{E \in B(C)}Z^E$ such that $\pi_{E}^{\ast} Z \cap E_{red}
\smallsetminus \pi_{E}^{\ast} C \cap E_{red} = \pi_{E}^{\ast} Z^{E}
\cap E_{red}$ by properties d2-d4.

Assume first that $Z$ is a perfect adjoint curve of $C$. Then
$\pi_C$ is a reduction of singularities of $Z \cup C$. Hence the
irreducible components of $Z$ are non-singular and its number is
equal to the multiplicity $m_0(Z)$. Moreover, the property d4.
implies that $\pi_E$ is a reduction of singularities of $Z^E$ and
the number of points of $\pi_{E}^{\ast}Z^{E} \cap E_{red}$ is equal
to $m_0(Z^{E}) =b_E-1$ since $Z^{E}$ only cuts $E_{red}$ by d3.

Reciprocally, assume that the set $\pi_E^* Z^{E} \cap E_{red}$ has
exactly $b_E-1$ points for each bifurcation divisor $E$ of $G(C)$.
This implies that $Z^E$ has $b_E-1$ irreducible components which are
non-singular and that $\pi_E$ is a reduction of singularities of
$Z^E$. Then, from the equalities $\pi_C^* Z^E \cap E = \pi_C^* Z
\cap E$ and $\pi_C^* Z^E \cap \pi_C^{\ast}C=\emptyset$, we deduce
that $\pi_C$ is a reduction of singularities of $Z \cup C$.
\end{proof}
The next corollary gives a characterization of a  perfect adjoint
curve of a given curve $C$  in terms of the equisingularity data of
$C$, when $C$ has only non-singular branches.
\begin{Corollary}\label{cor:perfadj-nosing}
Consider a curve  $C$  with only non-singular branches and let $Z
=\cup_{E \in B(C)}Z^{E}$ be the decomposition of a strict adjoint
curve $Z$ of $C$. The curve $Z$ is perfect adjoint curve of $C$  if
and only if each curve $Z^{E}$ is composed by $b_E-1$ irreducible
components $\{ \zeta_{i}^{E}\}_{i=1}^{b_E-1}$ with ${\mathcal
C}(\zeta_{i}^{E},\zeta_{j}^E)=v(E)$ for $i \neq j$.
\end{Corollary}

In particular, the corollary above implies that $G(C \cup Z)$ is
obtained from $G(C)$ by adding $b_E-1$ arrows to each bifurcation
divisor $E$ of $G(C)$ and this property characterizes the fact of
$Z$ being a perfect adjoint of $C$, when $C$ has only non-singular
branches. Hence, it is clear that $\epsilon(C \cup Z)$ does not
depend on $Z$ and we denote $\chi_C=\epsilon(C \cup Z)$.

In the general case of a curve $C$ with singular branches, it is not
possible to define $\chi_C$ in a compatible way with
$C$-ramifications. Since this situation needs a more detailed
treatment, we shall consider it in section~\ref{sec:tipo-amable}.

\section{Local invariants and polar curves}\label{sec:loc-inv}
Let $\mathbb F$ be the space of singular foliations of $({\mathbb
C}^2,0)$, that is, an element ${\mathcal F} \in {\mathbb F}$ is
defined by a 1-form $\omega=0$, with $\omega = A dx + B dy$,  $A, B
\in {\mathbb C}\{x,y\}$ and $A(0)=B(0)=0$. Given a plane curve $C
\subset ({\mathbb C}^2,0)$, we denote by ${\mathbb F}_C$ the
sub-space of ${\mathbb F}$ composed by the foliations which have $C$
as a curve of separatrices.

For a direction $[a:b] \in {\mathbb P}_{\mathbb C}^1$, the polar
curve $\Gamma({\mathcal F};[a:b])$ is the curve
$$\Gamma=\{ a A (x,y) + b B(x,y) =0\}.$$
We denote by $\Gamma_{\mathcal F}$ a generic polar when the
direction $[a:b]$ is not needed. Then the multiplicity
$m_0(\Gamma_{\mathcal F})$ of $\Gamma_{\mathcal F}$ at the origin
coincides with the multiplicity $\nu_0({\mathcal F})$ of $\mathcal
F$ at the origin. Recall that, if $\mathbb G$ is the space of
generalized curve foliations of $({\mathbb C}^2,0)$ and ${\mathbb
G}_C = {\mathbb F}_C \cap {\mathbb G}$, we have that
$\nu_0({\mathcal F})=m_0(C)-1$ for any ${\mathcal F} \in {\mathbb
G}_C$.

The Newton polygon ${\mathcal N}({\mathcal F};x,y) ={\mathcal
N}(\omega;x,y)$ of $\mathcal F$ is defined as the one of the ideal
generated by $x A $ and $y B$. More precisely, if we write $\omega =
\sum_{i,j} \omega_{ij}$ with
\begin{equation}\label{eq:1}
\omega_{ij} = A_{ij} x^{i-1} y^{j} dx + B_{ij} x^{i} y^{j-1} dy,
\end{equation}
and we put $\Delta(\omega)=\{ (i,j) \ : \ \omega_{ij}\neq 0\}$, then
${\mathcal N}({\mathcal F};x,y)$ is the convex envelop of
$\Delta(\omega) + {\mathbb R}_{\geq 0}^2$. In the case of an
analytic function $f = \sum_{ij} f_{ij}x^i y^j$, we define
$\Delta(f) = \{ (i,j) \ : \ f_{ij} \neq 0\}$ and then the Newton
polygon ${\mathcal N}(C;x,y)$ of the curve $C=(f=0)$ is the convex
envelop of $\Delta(f) + {\mathbb R}_{\geq 0}^2$. In particular, if
${\mathcal F} \in {\mathbb G}_C$, then ${\mathcal N}({\mathcal
F};x,y)$ coincides with ${\mathcal N}(C;x,y)={\mathcal N}(df;x,y)$.

From now on we will always assume that we chose coordinates $(x,y)$
such that $x=0$ is not tangent to the curve $C$ of separatrices. In
particular this implies that the first side of the Newton polygon
${\mathcal N}({\mathcal F};x,y)$ has slope greater or equal to $-1$.

Let us recall the relationship between Newton polygon and infinitely
near points of a curve since it will be useful  in the sequel. First
we introduce some notations
\begin{Notation}
Let $C$ be a curve with only non-singular branches and  $\pi_{C} : M
\rightarrow ({\mathbb C}^2,0)$ be the minimal reduction of
singularities of $C$. Given an irreducible component $E$ of
$\pi_C^{-1}(0)$ with $v(E)=p$, the morphism $\pi_{E} : M_E
\rightarrow ({\mathbb C}^{2},0)$ is a composition of $p$ blowing-ups
of points
$$
({\mathbb C}^2,0) \stackrel{\sigma_1}{\longleftarrow} (X_1,P_1)
\leftarrow \cdots \leftarrow (X_{p-1},P_{p-1})
\stackrel{\sigma_p}{\longleftarrow} X_p=M_E.
$$
If $(x,y)$ are coordinates in $({\mathbb C}^2,0)$ there is a change
of coordinates $(x,y)=(\tilde{x}, \tilde{y}+
\varepsilon(\tilde{x}))$, with $\varepsilon(x)= a_1 x + \cdots +
a_{p-1} x^{p-1}$,  such that the blowing up $\sigma_j$ is given by
$x_{j-1}=x_j$, $y_{j-1}=x_j y_j$, for $j=1,2, \ldots, p$, where
$(x_j,y_j)$ are coordinates centered at $P_j$ and
$(x_0,y_0)=(\tilde{x},\tilde{y})$. We say that
$(\tilde{x},\tilde{y})$ are {\em coordinates in\/} $({\mathbb
C}^2,0)$ {\em adapted to\/} $E$.
\end{Notation}

Consider now a plane curve $\gamma \subset ({\mathbb C}^2,0)$ with
only non-singular irreducible components and let $\pi_{\gamma}: X
\rightarrow ({\mathbb C}^2,0)$ be its minimal reduction of
singularities. Take $E$ an irreducible component of
$\pi_{\gamma}^{-1}(0)$ with $v(E)=p$ and choose $(x,y)$ coordinates
adapted to $E$. Assume that $\gamma = (f(x,y)=0)$ with $f(x,y)=
\sum_{i,j} f_{ij} x^{i} y^{j} \in {\mathbb C}\{x,y\}$. Since $(x,y)$
are adapted to $E$, then  there exists a side $L$ of ${\mathcal
N}(\gamma;x,y)$ with slope $-1/p$. Let $i+pj=k$ be the line which
contains $L$ and put
$$
In_{p}(f;x,y) = \sum_{i+pj=k} f_{ij} x^i y^j.
$$
Take now $(x_p,y_p)$ coordinates in the first chart of $E_{red}$
with $\pi_{E}(x_p,y_p)=(x_{p},x_{p}^{p} y_{p})$ and
$E_{red}=(x_p=0)$. Thus, a simple calculation shows that the points
of $\pi_{E}^{\ast} \gamma \cap E_{red}$ are given by $x_p=0$ and
$\sum_{i+pj=k} f_{ij} y^{j}=0$. We conclude that the points of
$\pi_{E}^* \gamma \cap E_{red}$ are determined by $In_p(f;x,y)$ and
reciprocally.

Consequently, the following  result which  describes the Newton
polygon of a generic polar curve $\Gamma_{\mathcal F}$ will be
useful to determine the infinitely near points of $\Gamma_{\mathcal
F}$.
\begin{Lemma}[\cite{Cor-03}]\label{lem:polNew:F-Gamma}
Consider a foliation ${\mathcal F} \in {\mathbb F}$ and let $L$ be a
side of ${\mathcal N}({\mathcal F};x,y)$ with slope $-1/\mu$ where
$\mu \in {\mathbb Q}$ and $\mu \geq 1$. If $i+\mu j =k$ is the
equation of the line which contains $L$, then
$$
{\mathcal N}(\Gamma_{\mathcal F};x,y) \subset \{ (i,j) \ : \ i + \mu
j \geq k-\mu\}.
$$
More precisely, if $\mu >1$ then $\Delta(B) \subset \{i+\mu j \geq
k-\mu\}$ and $\Delta(A) \subset \{i+\mu j > k-\mu\}$.
\end{Lemma}
However the result above does not provide enough information to
obtain a description of the equisingularity type  of
$\Gamma_{\mathcal F}$. If we want to control the slopes of
${\mathcal N}(\Gamma_{\mathcal F};x,y)$ we need to know the
``contribution" in the points of the sides of ${\mathcal
N}({\mathcal F};x,y)$. Recall that  a point $(i,j) \in
\Delta(\omega)$ is said to be a {\em contribution\/} of $B$ if
$B_{ij} \neq 0$ in the expression \eqref{eq:1}, i.e., if $(i,j) \in
\Delta(yB)$.

Thus to get a more precise description of the Newton polygon
${\mathcal N}(\Gamma_{\mathcal F};x,y)$ we need to consider
foliations in ${\mathbb G}_{C}^*$ since the contributions on the
sides of the Newton polygon of a foliation have a direct
relationship with the values of the Camacho-Sad indices at the
infinitely near points of $\mathcal F$ as it is explained in the
next proposition.

Recall that, if $S=(y=0)$ is a non-singular separatrix of $\mathcal
F$, then the {\em Camacho-Sad index of $\mathcal F$ relative to $S$
at the origin\/} is given by
\begin{equation}\label{eq:calculo-indice}
I_0 ({\mathcal F},S)= -\text{Res}_0 \frac{a(x,0)}{B(x,0)}
\end{equation}
where the $1$-form $\omega$ defining $\mathcal F$ is written as
$\omega = y a(x,y) dx + b(x,y) dy$ (see \cite{Cam-S}). Then we have
the following result:

\begin{Proposition}[\cite{Cor-03}]\label{indice:-1}
Consider a foliation  ${\mathcal F} \in {\mathbb G}$ and take  a
side $L$ of ${\mathcal N}({\mathcal F})$ with slope $-1/p$, $p \in
{\mathbb N}$. If $L$ has no contribution of $B$ in its highest
vertex, then there is a corner in the reduction of singularities of
$\mathcal F$ with Camacho-Sad index equal to $-1$.
\end{Proposition}
In particular, given  a foliation ${\mathcal F} \in {\mathbb
G}_{C}^{\ast}$ such that the curve $C$ has only non-singular
irreducible components, the result above implies that
\begin{quote}
if ${\mathcal N}({\mathcal F};x,y)$ has $s$ sides $L_j$ with slopes
$-1/p_j$, $p_j \in {\mathbb N}$, $j=1, \ldots ,s$ and $p_1 < p_2 <
\cdots < p_s$, then the first $s-1$ sides of ${\mathcal
N}(\Gamma_{\mathcal F};x,y)$ are obtained from the ones of
${\mathcal N}({\mathcal F};x,y)$ by a vertical translation of one
unit and the other ones have slope  $\geq -1/p_s$.
\end{quote}
These results describing the Newton polygon of $\Gamma_{\mathcal F}$
are key in the proof of the decomposition theorem:
\begin{Theorem}[of decomposition \cite{Cor-03}]\label{th-decomposition}
Consider a foliation ${\mathcal F} \in {\mathbb G}_{C}^{\ast}$ and
$\Gamma_{\mathcal F}$ a generic polar curve of $\mathcal F$. Given
any $C$-ramification $\rho: ({\mathbb C}^2,0) \rightarrow ({\mathbb
C}^2,0)$, the curve $\rho^{-1}\Gamma_{\mathcal F}$ is a strict
adjoint of $\rho^{-1}C$.
\end{Theorem}
By the results in section~\ref{sec:adjuntos}, we deduce that there
is a unique decomposition $\rho^{-1} \Gamma_{\mathcal F} =
\cup_{\tilde{E} \in B(\tilde{C})} \tilde{\Gamma}^{\tilde{E}}$, with
$\tilde{C}=\rho^{-1}C$, satisfying the properties d1-d4, (D-i) and
(D-ii). Moreover, the curve $\Gamma_{\mathcal F}$ can also be
decomposed in unique way as
$$
\Gamma_{\mathcal F} = \bigcup_{E \in B(C)} \Gamma^E
$$
satisfying properties D1-D5, (D-i) and (D-ii) in
section~\ref{sec:adjuntos}.

Observe now that the property of being a strict adjoint  of a curve
$C$ does not determine the equisingularity type of the adjoint
curve: for instance, if  $C$ is the union of 3 lines, then there are
infinite many  possible equi\-singularity types for its strict
adjoint curves. However, the number of possible equisingularity
types is finite when considering polar curves.

\begin{Theorem}\label{th:1}
There exists a finite number of equisingularity types $\wp({\mathcal
F})$ for a foliation ${\mathcal F} \in {\mathbb G}_{C'}^{\ast}$ and
 any curve $C'$ with $\epsilon(C')=\epsilon(C)$.
\end{Theorem}
\begin{proof}
Let $\mathcal F$ be a foliation in ${\mathbb G}_C^*$ and consider a
generic polar curve $\Gamma=\Gamma_{\mathcal F}$ of $\mathcal F$. It
is clear that the number of irreducible components of $\Gamma$ is
lower than or equal to the multiplicity $m_0(\Gamma)=m_0(C)-1$.

Consider a ramification $\rho: ({\mathbb C}^2,0) \rightarrow
({\mathbb C}^2,0)$ transversal to $C$ and such that $\rho^{-1} C$
and $\rho^{-1}\Gamma$ have non-singular irreducible components.
Let us prove that given any two irreducible components $\sigma,
\sigma'$ of $\rho^{-1}\Gamma$ the coincidence ${\mathcal
C}(\sigma,\sigma')$ is bounded in terms of the equisingularity data
of $\rho^{-1}C$. In particular, this implies that there is only a
finite number of possibilities for the characteristic exponents of
the branches of $\Gamma$ and for the coincidence between two
branches of $\Gamma$ once the equisingularity type of $C$ is fixed
(see appendix \ref{ap:ramificacion}). Hence, the number of possible
equisingularity types for $\Gamma$ is finite. Moreover, since the
coincidences between the irreducible components of $\Gamma$ and $C$
are determined by $\epsilon(C)$, the result follows straightforward.

Let $p=\sup_{\sigma, \sigma'} {\mathcal C}(\sigma,\sigma')$ where
$\sigma, \sigma'$ vary within the irreducible components of
$\rho^{-1}\Gamma$; observe that $p \in {\mathbb N}$. If $p \leq
\sup_{\alpha, \alpha'} {\mathcal C}(\alpha, \alpha')$ for $\alpha,
\alpha'$  among the irreducible components of $\rho^{-1}C$ we
finish. Otherwise let $\sigma_0,\sigma_0'$ be two irreducible
components of $\rho^{-1}\Gamma$ such that ${\mathcal
C}(\sigma_0,\sigma_0')=p$. In particular, by property (D-ii) of the
decomposition of $\rho^{-1}\Gamma$, we have that
$\mu=\sup_{\alpha}{\mathcal
C}(\sigma_0,\alpha)=\sup_{\alpha}{\mathcal C}(\sigma_0',\alpha) < p$
where $\alpha$ varies within the irreducible components of
$\rho^{-1}C$.

Take $(x,y)$ coordinates in $({\mathbb C}^2,0)$ such that the
coincidence of the axis $y=0$ with the curves $\sigma_0$ and
$\sigma_0'$ is equal to $p$. This implies that the last side
$L_{\tilde{\Gamma}}$ of the Newton polygon ${\mathcal
N}(\rho^{-1}\Gamma;x,y)$ has a slope equal to $-1/p$. Moreover, the
last side $L_{\mathcal F}$ of ${\mathcal N}(\rho^{\ast}{\mathcal
F};x,y)$ has a slope equal to $-1/\mu$.

Let $i + \mu j = k$ be the line which contains $L_{\mathcal F}$ and
$(l_1,h_1)$  be the highest  vertex of $L_{\mathcal F}$ (note that
$h_1 \geq 3$). The previous results concerning the behaviour of the
Newton polygon ${\mathcal N}(\rho^{-1}\Gamma;x,y)$ imply that a
point $(i,j)$ on $L_{\tilde{\Gamma}}$ must verify the following
conditions
\begin{equation*}
\left\{ \begin{array}{ll} \ 0 \leq j \leq h_1-1  & \text{ by prop.
\ref{indice:-1}}; \\
i + \mu j \geq k -\mu & \text{ by lemma \ref{lem:polNew:F-Gamma}}; \\
i + \frac{k-l_1-1}{h_1-1} j \leq k-1  & \text{ since $(l_1,h_1-1),
(k-1,0) \in \Delta(\rho^{-1}\Gamma$)}.
\end{array} \right.
\end{equation*}
Thus there exists only a finite number of possible values for $p$.
Moreover, from the inequalities above we deduce that $\mu \leq p <2
\mu$. The next picture illustrate the situation: the side
$L_{\tilde{\Gamma}}$ must be contained in the grey region with slope
equal to $-1/p$, $p \in {\mathbb N}$.
\begin{center}
\begin{texdraw}
\drawdim mm \setgray 0 \linewd 0.1
\def\ticklab(#1 #2){\move(#1 #2)
  \bsegment
  \lvec(0 1)
  \esegment}
\def\vticklab(#1 #2){\move(#1 #2)
  \bsegment
  \lvec(1 0)
  \esegment}
\arrowheadtype t:F \arrowheadsize l:2 w:1

\move (10 10) \avec (10 40) \move(7 38) \htext{\tiny{$j$}} \move(10
10) \avec (100 10) \move(98 7) \htext{\tiny{$i$}} \move(90 7)
\htext{\tiny{$k$}}

\linewd 0.3 \move(20 35)
\rlvec (10 -10)
\move (30 25) \fcir f:0 r:0.5 \rlvec (60 -15) \move (90 10) \fcir
f:0 r:0.5

\move (25 33) \htext{\tiny{${\mathcal N}(\rho^{\ast}{\mathcal F})$}}
\move (30 26) \htext{\tiny{$(l_1,h_1)$}}

\linewd 0.1 \ticklab(15 9.5) \ticklab(20 9.5) \ticklab(25 9.5)
\ticklab(30 9.5) \ticklab(35 9.5) \ticklab(40 9.5) \ticklab(45 9.5)
\ticklab(50 9.5) \ticklab(55 9.5) \ticklab(60 9.5) \ticklab(65 9.5)
\ticklab(70 9.5) \ticklab(75 9.5) \ticklab(80 9.5) \ticklab(85 9.5)
\ticklab(90 9.5) \ticklab(95 9.5)

\vticklab(9.5 15) \vticklab(9.5 20) \vticklab(9.5 25) \vticklab(9.5
30) \vticklab(9.5 35)

\setgray 0.5 \linewd 0.4
 \move(20 30)
\rlvec (10 -10)
 \lpatt(1 1) \rlvec (40 -10) \lpatt()
\move (85 10) \fcir f:0.5 r:0.5 \move (30 20) \lpatt(1 1) \rlvec (55
-10) \lpatt() \move (30 20) \fcir f:0.4 r:0.5

\move(30 20) \rlvec(40 -10) \rlvec(15 0) \rlvec(-55 10) \ifill f:0.7

\move (60 18) \htext{\tiny{$L_{\mathcal F}$}}
 \move (13 23) \htext{\tiny{${\mathcal N}(\rho^{-1}\Gamma)$}}
\end{texdraw}
\end{center}
\end{proof}
Among all the possible equisingularity types  $\wp({\mathcal
F})=\epsilon(\Gamma_{\mathcal F} \cup C)$ for a fixed
equisingularity type $\epsilon(C)$, there is one which can be
considered as the ``minimal" one satisfying the decomposition
theorem. Next sections will be devoted to characterize foliations
such that $\wp({\mathcal F})$ is the minimal one.

\section{Non-singular branches}\label{sec:non-sp}
In this section we consider a curve $C = \cup_{i=1}^r C_i$ with only
non-singular irreducible components and we study under what
conditions a generic polar curve $\Gamma_{\mathcal F}$ of a
foliation ${\mathcal F} \in {\mathbb G}_{C}^*$ is a perfect adjoint
of $C$. Denote by ${\mathbb G}_{C,\lambda}^*$ the space of
foliations ${\mathcal F} \in {\mathbb G}_{C}^*$ such that
$\lambda({\mathcal F})=\lambda$. Let $U_C \subset {\mathbb
P}^{r-1}_{\mathbb C}$ be the set defined by
\begin{quote}
$\lambda \in U_C$ if there exists ${\mathcal F} \in {\mathbb
G}_{C,\lambda}^{\ast}$ with $\wp({\mathcal F})=\chi_C$.
\end{quote}
Then we have
\begin{Theorem}[of genericity]\label{th:U_C}
The set $U_C$ is a non-empty Zariski open set. Moreover, for each
${\mathcal F} \in {\mathbb G}_{C,\lambda}^{\ast}$ with $\lambda \in
U_C$ we have that $\wp({\mathcal F})=\chi_C$.
\end{Theorem}
\begin{Definition}
A foliation ${\mathcal F} \in {\mathbb G}_{C}^*$ is  {\em
Zariski-general\/} if $\lambda({\mathcal F}) \in U_C$.
\end{Definition}

Denote by ${\mathcal L}_{\lambda}$ a logarithmic foliation in
${\mathbb G}_C$ with $\lambda= (\lambda_1, \ldots,\lambda_r) \in
{\mathbb P}_{\mathbb C}^{r-1}$. We define the set
$$
U_{C}^{log} = \{ \lambda \in {\mathbb P}_{\mathbb C}^{r-1} \ : \
{\mathcal L}_{\lambda} \in {\mathbb G}_{C}^{\ast} \text{ and }
\wp({\mathcal L}_{\lambda}) = \chi_C \}.
$$
It is clear that $U_{C}^{\log} \subset U_C$. Let us prove the
following result
\begin{Proposition}\label{prop:U_C-log}
The set $U_C^{log}$ is a non-empty Zariski open set of ${\mathbb
P}_{\mathbb C}^{r-1}$.
\end{Proposition}
\begin{proof}
We note first that the equisingularity type of a generic polar curve
of a logarithmic foliation ${\mathcal L}_{\lambda} \in {\mathbb
F}_{C}$ does not depend on the equations of $C=\cup_{i=1}^{r}C_i$
chosen to define ${\mathcal L}_{\lambda}$ (see prop. 3.8 of
\cite{Cor-03}). So we can assume that  ${\mathcal L}_{\lambda}$ is
defined  by $\omega_{\lambda}=0$ with
\begin{equation}\label{eq:fol-log-nosing}
\omega_{\lambda} = \prod_{i=1}^{r} (y-\eta_i(x)) \sum_{i=1}^{r}
\lambda_i \frac{d(y-\eta_{i}(x))}{(y-\eta_i(x))},
\end{equation}
where the curve $C_i$ is defined by $(y-\eta_i(x)=0)$ and $\eta_i(x)
= \sum_{j=1}^{\infty} a_{j}^{i} x^{j} \in {\mathbb C}\{x\}$.
Moreover, for a direction $[a:b] \in {\mathbb P}_{\mathbb C}^1$, the
polar curve $\Gamma({\mathcal L}_{\lambda};[a:b])$ is given by
\begin{equation}\label{eq:polar-log-nosing}
\sum_{i=1}^{r} \lambda_i \prod_{j \neq i} (y - \eta_j(x)) (-a
\eta_{i}'(x) + b)=0
\end{equation}
and we denote by $\Gamma^{\lambda}_{[a:b]}$ a generic polar curve of
${\mathcal L}_{\lambda}$.

The first condition over $\lambda$ to belong to $U_C^{log}$ is that
${\mathcal L}_{\lambda} \in {\mathbb G}_C^*$ but this is equivalent
to $\sum_{i=1}^{r} k_i \lambda_i \neq 0$ where $k \in
R_{\epsilon(C)}$ and $R_{\epsilon(C)}$ is a finite set of
reso\-nances (see \cite{Cor-03} for a detailed description of
$R_{\epsilon(C)}$). Now, for each bifurcation divisor $E$ of $G(C)$,
we define $U_C^E$ to be the set of $\lambda \in {\mathbb P}_{\mathbb
C}^{r-1}$ such that $\pi_{E}^{\ast} \Gamma_{[a:b]}^{\lambda} \cap
E_{red} \smallsetminus \pi_{E}^* C \cap E_{red}$ has exactly $b_E-1$
different points, and we will prove that $U_C^E$ is a non-empty
Zariski open set. Using the criterion given in
proposition~\ref{criterio-adjunto-perfecto}, we obtain that
$U_{C}^{log}$ is equal to
$$
U_C^{log}=\{ \lambda \in {\mathbb P}_{\mathbb C}^{r-1} \ : \ \lambda
\in \bigcap_{E \in B(C)} U_C^E \text{ and } \sum_{i=1}^r k_i
\lambda_i \neq 0 \text{ for } k \in R_{\epsilon(C)} \}
$$
which is a non-empty Zariski open set.

Take a bifurcation divisor $E$ of $G(C)$ with $v(E)=p$ and let us
prove that each $U_C^E$ is a non-empty Zariski open set. Let
$\pi_{E}: M_E \rightarrow ({\mathbb C}^2,0)$ be the reduction of
$\pi_C$ to $E$. Since the equisingularity type of a generic polar
curve of a foliation does not depend on the coordinates (see
\cite{Cor-03}, \S 2), we can assume that the coordinates $(x,y)$ are
adapted to $E$. Take $(x_p,y_p)$ coordinates in the first chart of
$E_{red}\subset M_E$ such that $\pi_{E}(x_p,y_p)=(x_p,x_p^{p} y_p)$
and $E_{red}=(x_p=0)$. If the strict transform
$\pi_{E}^{\ast}{\mathcal L}_{\lambda}$ of ${\mathcal L}_{\lambda}$
is defined by $\omega_{\lambda}^E =0$ with
$$
\omega_{\lambda}^{E}=A_{\lambda}^E(x_p,y_p)dx_p + x_{p}
B_{\lambda}^{E} (x_p,y_p) dy_p,
$$
then the singular points of $\pi_{E}^{\ast}{\mathcal L}_{\lambda}$
in the first chart of $E_{red}$ are given by $x_p=0$ and
$A_{\lambda}^E(0,y_p)=0$. Let us compute the polynomials
$A_{\lambda}^{E}(0,y)$ and $B_{\lambda}^{E}(0,y)$.

We consider two situations:  $E$ being the first bifurcation divisor
of $G(C)$ or not. If $E$ is the first bifurcation divisor, then $E$
belongs to the geodesic of all the irreducible components of $C$.
Let $\{R_1^{E}, \ldots, R_{b_E}^{E}\}$ be the singular points of
$\pi_{E}^{\ast} {\mathcal L}_{\lambda}$ in the first chart of
$E_{red}$ where $R_i^{E}=(0,c_i^{E})$ in the coordinates
$(x_p,y_p)$.

Compute the strict transform of $\omega_{\lambda}$ by $\pi_E$ using
the expression in~\eqref{eq:fol-log-nosing} and the fact that
$\{R_{1}^{E}, \ldots, R_{b_E}^{E}\} = \pi_{E}^* C \cap E_{red}$,
thus we get that
\begin{align*}
A_{\lambda}^E(0,y) &= p \cdot (\sum_{i=1}^r \lambda_i)
\prod_{i=1}^{r}(y-a_{p}^{i}) = p \ (\sum_{i=1}^{r} \lambda_i)
\prod_{l=1}^{b_E} (y -c_{l}^{E})^{r_l} \\
B_{\lambda}^{E}(0,y) & = \sum_{i=1}^{r} \lambda_i \prod_{j \neq i}
(y-a_{p}^{j})
\end{align*}
where $r_l= m_{R_l^{E}}(\pi_{E}^{\ast}C)$; note that also $r_l=
\sharp\{j \ : \pi_{E}^{\ast}C_j \cap E_{red} = \{R_l^E\} \}$.

Let us now compute the strict transform of
$\Gamma_{[a:b]}^{\lambda}$ by $\pi_E$. By  the equation of
$\Gamma_{[a:b]}^{\lambda}$ given in~\eqref{eq:polar-log-nosing} and
lemma~\ref{lem:polNew:F-Gamma}, we obtain that the points of the set
$\pi_{E}^{\ast} \Gamma_{[a:b]}^{\lambda} \cap E_{red}$ are given by
$x_p=0$ and
\begin{equation}\label{eq:points-polar}
\left\{
  \begin{array}{ll}
    B_{\lambda}^{E}(0,y_p)=0, & \hbox{ if $p>1$;} \\
    a A_{\lambda}^{r-1}(1,y_p)+ b B_{\lambda}^{r-1}(1,y_p)=0, & \hbox{ if $p=1$,}
  \end{array}
\right.
\end{equation}
where $A_{\lambda}^{r-1}(x,y) dx + B_{\lambda}^{r-1}(x,y) dy$ is the
jet of order $\nu_0({\mathcal L}_{\lambda})=r-1$ of
$\omega_{\lambda}$. Hence we shall consider the two cases: $p>1$ and
$p=1$ to describe the set $\pi_{E}^{\ast} \Gamma_{[a:b]}^{\lambda}
\cap E_{red} \smallsetminus \pi_{E}^*C \cap E_{red}$.

By theorem~\ref{th-decomposition}, we know that
$m_{R_{i}^{E}}(\pi_{E}^{\ast} \Gamma_{[a:b]}^{\lambda})=r_l-1$ and
consequently,   the polynomial $\prod_{l=1}^{b_E}(y-c_l^E)^{r_l-1}$
divides the polynomials in \eqref{eq:points-polar}. In particular,
the points of $\pi_{E}^{\ast} \Gamma_{[a:b]}^{\lambda} \cap E_{red}
\smallsetminus \pi_{E}^*C \cap E_{red}$ are given by $x_p=0$ and
$H_{\lambda}^E (y_p)=0$ with
\begin{equation*}
H_{\lambda}^E(y)= \left\{
  \begin{array}{ll}
    B_{\lambda}^{E}(0,y)/\prod_{l=1}^{b_E}(y-c_{l}^{E})^{r_l-1}, & \hbox{if $p>1$;} \\
    (a A_{\lambda}^{r-1}(1,y)+ b B_{\lambda}^{r-1}(1,y))/\prod_{l=1}^{b_E}(y-c_{l}^{E})^{r_l-1},
& \hbox{if $p=1$.}
  \end{array}
\right.
\end{equation*}
The degree of $H_{\lambda}^E(y)$ as a polynomial in $y$ is equal to
$b_E-1$ and its coefficients depend linearly on the $\lambda_i$. Let
us study   the two cases $p>1$ and $p=1$.

{\sc Case $p>1$:} Let $D^E(\lambda)$ be the discriminant of
$H_{\lambda}^E(y)$ as a polynomial in $y$. Thus, the polynomial
$H_{\lambda}^E(y)$ has $b_E-1$ different roots if and only if
$D^E(\lambda) \neq 0$. Note that $D^E(\lambda) \not \equiv 0$ since
$D^E (1,0,\ldots,0) \neq 0$. Thus, the set $U_C^E$ is equal to  the
non-empty Zariski open set ${\mathbb P}_{\mathbb C}^{r-1}
\smallsetminus \{D^E =0\}$.

{\sc Case $p=1$:} The exceptional divisor $E$ coincides with $E_1$
and the coordinates $(x,y)$ are adapted to $E_1$.
From~\eqref{eq:fol-log-nosing} we get that
\begin{align*}
A_{\lambda}^{r-1}(1,y) &=-\sum_{i=1}^{r} \lambda_i a_{1}^{i}
\prod_{j \neq i} (y -a_{1}^{j}) \\
B_{\lambda}^{r-1} (1,y) & = B_{\lambda}^{E_1}(0,y) = \sum_{i=1}^r
\lambda_i \prod_{j \neq i} (y-a_1^j).
\end{align*}
Thus the polynomial $H_{\lambda}^{E_1}(y)$ can be written as follows
$$
H_{\lambda}^{E_1}(y) = \frac{a A_{\lambda}^{r-1} (1,y) + b
B_{\lambda}^{r-1}(1,y)}{\prod_{l=1}^{b}(y-c_l^{E_1})^{r_l-1}} = a
A_{\lambda}^\natural(y) + b B_{\lambda}^{\natural}(y).
$$
Let us show that $H_{\lambda}^{E_1}(y)$ has $b_{E_1}-1$ different
roots. It is clear that
$$
A_{\lambda}^{E_1}(0,y) = A_{\lambda}^{r-1}(1,y) + y
B_{\lambda}^{r-1}(1,y) = (\sum_{i=1}^{r} \lambda_i)
\prod_{l=1}^{b_{E_1}} (y - c_l^{E_1})^{r_l}
$$
and then $A_{\lambda}^{\natural}(y) + y B_{\lambda}^{\natural}(y) =
(\sum_{i=1}^r \lambda_i) \prod_{l=1}^{b_{E_1}}(y-c_l^{E_1})$. In
particular, we deduce that $A_{\lambda}^{\natural}(y)$ and
$B_{\lambda}^{\natural}(y)$ do not have common roots. In fact, the
only possible common roots are the elements of the set
$\{c_l^{E_1}\}_{l=1}^{b_{E_1}}$, but if $c_l^{E_1}$ is a common root
of both polynomials then it is also a root of $H_{\lambda}^{E_1}(y)$
in contradiction with theorem~\ref{th-decomposition}. Thus for $a,b$
generic, the polynomial $H_{\lambda}^{E_1}(y)$ has $b_{E_1}-1$
different roots and hence $U_C^{E_1}={\mathbb P}_{\mathbb C}^{r-1}$.

We consider now the case of $E$ being any bifurcation divisor. Put
$I=\{1,2,\ldots, r\}$ and $I^E=\{i \in I \ : \ E \text{ belongs to
the geodesic of } C_i\}$. We can write $\omega_{\lambda}=
\omega_{\lambda}^*+\omega_{\lambda}^{**}$ where
\begin{align*}
\omega_{\lambda}^* &= \prod_{i \in I^E} (y-\eta_i(x)) \sum_{j \in I
\smallsetminus I^E} \lambda_j \prod_{l \in I \smallsetminus I^E
\atop l  \neq j} (y-\eta_l(x)) (-\eta_j'(x) dx + dy), \\
\omega_{\lambda}^{**} &= \prod_{i \in I \smallsetminus I^E}
(y-\eta_i(x)) \sum_{j \in I^E} \lambda_j \prod_{l \in  I^E \atop l
\neq j} (y-\eta_l(x)) (-\eta_j'(x) dx + dy).
\end{align*}
If we compute the strict transform $\omega_{\lambda}^E$ of
$\omega_{\lambda}$ by $\pi_E$,  we get that the polynomials
$A_{\lambda}^{E}(0,y)$ and $B_{\lambda}^{E}(0,y)$ are given by
\begin{equation*}
A_{\lambda}^{E}(0,y) = C \cdot \prod_{i \in I^{E}} (y-a_p^{i});  \ \
B_{\lambda}^{E}(0,y) = C' \cdot \sum_{i \in I^{E}} \prod_{i \in
I^{E} \atop j \neq i} (y- a_p^j)
\end{equation*}
where $C, C'$ are non-zero constants. Thus the  set $U_{C}^{E}$ is
defined in a similar way to the case of $E$ being the first
bifurcation divisor with $p>1$.

We conclude that $U_{C}^{log}$ is a non-empty Zariski open set
because it is a finite intersection of non-empty Zariski open sets.
\end{proof}
The next lemma concerns the infinitely near points of  generic polar
curves and, in particular, it allows to show the equality of the
sets $U_{C}$ and $U_{C}^{log}$.
\begin{Lemma}\label{lema-puntos-polarF-L}
Consider two foliations ${\mathcal F}, {\mathcal L}_{\lambda} \in
{\mathbb G}_{C,\lambda}^*$. Let $\Gamma_{[a:b]}^{\mathcal F}$ and
$\Gamma_{[a:b]}^{{\mathcal L}_{\lambda}}$ be  generic polar curves
of $\mathcal F$ and ${\mathcal L}_{\lambda}$ respectively. Then, for
each irreducible component $E$ of $\pi_{C}^{-1}(0)$, we have that
$$
\pi_{E}^{\ast} \Gamma_{[a:b]}^{\mathcal F} \cap E_{red} =
\pi_{E}^{\ast} \Gamma_{[a:b]}^{{\mathcal L}_{\lambda}} \cap E_{red}
$$
and the multiplicities satisfy that
$m_P(\pi_{E}^*\Gamma_{[a:b]}^{\mathcal F})=m_P(\pi_{E}^*
\Gamma_{[a:b]}^{{\mathcal L}_{\lambda}})$ at each point $P \in
\pi_{E}^{\ast} \Gamma_{[a:b]}^{\mathcal F} \cap E_{red}$. Moreover,
if $E \neq E_1$, the sets above does not depend on $[a:b]$, that is,
$$
\pi_{E}^{\ast} \Gamma_{[a:b]}^{\mathcal F} \cap E_{red} =
\pi_{E}^{\ast} \Gamma_{[a':b']}^{\mathcal F} \cap E_{red} =
\pi_{E}^{\ast} \Gamma_{[a:b]}^{{\mathcal L}_{\lambda}} \cap E_{red}
=\pi_{E}^{\ast} \Gamma_{[a':b']}^{{\mathcal L}_{\lambda}} \cap
E_{red}
$$
for all $[a:b], [a':b']$ generic.
\end{Lemma}
\begin{proof}
Take an irreducible component $E$ of $\pi_{C}^{-1}(0)$ and let
$\pi_E : M_E \rightarrow ({\mathbb C}^2,0)$ be the reduction of
$\pi_C$ to $E$. If $E$ is not a bifurcation divisor, then $\pi_{E}^*
\Gamma_{[a:b]}^{\mathcal F} \cap E_{red}$ and $\pi_{E}^*
\Gamma_{[a:b]}^{{\mathcal L}_{\lambda}} \cap E_{red}$ coincide with
$\pi_{E}^* C \cap E_{red}$ because $\Gamma_{[a:b]}^{\mathcal F}$ and
$\Gamma_{[a:b]}^{{\mathcal L}_{\lambda}}$ are strict adjoint curves
of $C$; in particular, the points of the set $\pi_{E}^{\ast}
\Gamma_{[a:b]}^{\mathcal F} \cap E_{red}$ does not depend on
$[a:b]$. Moreover, $m_P(\pi_{E}^* \Gamma_{[a:b]}^{\mathcal F})
=m_P(\pi_E^* \Gamma_{[a:b]}^{{\mathcal L}_{\lambda}}) =
m_P(\pi_{E}^* C)-1$ at each point $P \in \pi_{E}^* C \cap E_{red}$
by theorem~\ref{th-decomposition}.

Assume now that $E$ is a bifurcation divisor with $v(E)=p$. In order
to simplify notations, we suppose that $E$ is the first bifurcation
divisor and that the coordinates $(x,y)$ are adapted to $E$;
otherwise we work in a similar way  as in the proof of
proposition~\ref{prop:U_C-log}. Consider two $1$-forms
$\omega_{\mathcal F}=A_{\mathcal F}(x,y) dx + B_{\mathcal F}(x,y)dy$
and $\omega_{\mathcal L}= A_{\mathcal L}(x,y) dx + B_{\mathcal
L}(x,y) dy$ such that $\mathcal F$ and ${\mathcal L}={\mathcal
L}_\lambda$ are defined by $\omega_{\mathcal F}=0$ and
$\omega_{\mathcal L}=0$ respectively.

Take $(x_p,y_p)$  coordinates in the first chart of $E_{red}$ such
that $\pi_{E}(x_p,y_p)=(x_p,x_p^p y_p)$ and $E_{red}=(x_p=0)$. Let
$\omega^{E}_{\mathcal F}$ and $\omega_{\mathcal L}^E$ be the strict
transforms  of $\omega_{\mathcal F}$ and $\omega_{\mathcal L}$ by
$\pi_E$ with
\begin{align}
\omega_{\mathcal F}^{E} & = A_{\mathcal F}^E (x_p,y_p) dx_p + x_p
B_{\mathcal F}^{E} (x_p,y_p) dy_p , \label{eq:omegaF-E} \\
\omega_{\mathcal L}^{E} & = A_{\mathcal L}^E (x_p,y_p) dx_p + x_p
B_{\mathcal L}^{E} (x_p,y_p) dy_p. \label{eq:omegaL-E}
\end{align}

Denote by $\{R_1^E,\ldots,R_{b_E}^E\}$ the points of the set
$\pi_{E}^*C \cap E_{red}$ and assume that each point
$R_{l}^{E}=(0,c_{l}^E)$ in the coordinates $(x_p,y_p)$. The singular
points of $\pi_{E}^*{\mathcal F}$ and $\pi_{E}^{\ast}{\mathcal L}$
in the first chart of $E_{red}$ coincide with the points of
$\pi_{E}^*C \cap E_{red}$ since ${\mathcal F}$ and ${\mathcal L}$
belong to ${\mathbb G}_{C}$. Moreover,
$m_{R_{i}^{E}}(\pi_{E}^*{\mathcal
F})=m_{R_{i}^{E}}(\pi_{E}^*{\mathcal L}) =
m_{R_{i}^{E}}(\pi_{E}^*C)$. Thus,  up to divide $\omega_{\mathcal
F}^{E}$ and $\omega_{\mathcal L}^E$ by a constant, we have that
\begin{equation}\label{eq:A-E}
A_{\mathcal F}^E(0,y)=A_{\mathcal L}^E(0,y) = \prod_{l=1}^{b_E}
(y-c_{l}^E)^{r_l}
\end{equation}
with $r_l=m_{R_{l}^{E}}(\pi_{E}^*C)$. By
theorem~\ref{th-decomposition}, we also have that
$m_{R_{l}^E}(\pi_{E}^* \Gamma_{[a:b]}^{\mathcal F}) =
m_{R_{l}^E}(\pi_{E}^* \Gamma_{[a:b]}^{\mathcal L})=
m_{R_{l}^E}(\pi_{E}^* C)-1$. Thus we only need to show that the sets
$\pi_{E}^* \Gamma_{[a:b]}^{\mathcal F} \cap E_{red} \smallsetminus
\pi_{E}^*C \cap E_{red}$ and $\pi_{E}^* \Gamma_{[a:b]}^{\mathcal L}
\cap E_{red} \smallsetminus \pi_{E}^*C \cap E_{red}$ coincide. Using
similar arguments as in the proof of proposition~\ref{prop:U_C-log},
we obtain that the points of $\pi_{E}^* \Gamma_{[a:b]}^{\mathcal F}
\cap E_{red} \smallsetminus \pi_{E}^*C \cap E_{red}$ are given by
$x_p=0$ and $H^E_{\mathcal F}(y_p)=0$ where
\begin{equation*}
H_{\mathcal F}^E(y)= \left\{
  \begin{array}{ll}
    B_{\mathcal F}^{E}(0,y)/\prod_{l=1}^{b_E}(y-c_{l}^{E})^{r_l-1}, & \hbox{if $p>1$;} \\
    (a A_{\mathcal F}^{r-1}(1,y)+ b B_{\mathcal F}^{r-1}(1,y))/\prod_{l=1}^{b_{E_1}}(y-c_{l}^{E_1})^{r_l-1},
& \hbox{if $p=1$,}
  \end{array}
\right.
\end{equation*}
and $A_{\mathcal F}^{r-1}(x,y) dx + B_{\mathcal F}^{r-1}(x,y) dy$ is
the jet of order $\nu_0({\mathcal F})=r-1$ of $\omega_{\mathcal F}$.
We obtain in a similar way a polynomial $H_{\mathcal L}^{E}(y)$ for
the foliation $\mathcal L$. In order to prove the lemma we only need
to show that the polynomials $H_{\mathcal F}^{E}(y)$ and
$H_{\mathcal L}^{E}(y)$ coincide.

Taking into account that $\mathcal L$ is a logarithmic model of
$\mathcal F$, we get that the Camacho-Sad indices
$I_{R_{l}^E}(\pi_E^* {\mathcal F},E_{red})$ and $I_{R_{l}^E}(\pi_E^*
{\mathcal L},E_{red})$ are equal for $l=1,\ldots,b_E$. From the
definition of the Camacho-Sad index given
in~\eqref{eq:calculo-indice} and equations~\eqref{eq:omegaF-E},
\eqref{eq:omegaL-E} we obtain that
\begin{align*}
I_{R_{l}^E}(\pi_E^* {\mathcal F},E_{red})  &=
\text{Res}_{y=c_{l}^{E}} \frac{-B_{\mathcal F}^{E}(0,y)}{A_{\mathcal
F}^{E}(0,y)}; \\
I_{R_{l}^E}(\pi_E^* {\mathcal L},E_{red}) &=
\text{Res}_{y=c_{l}^{E}} \frac{-B_{\mathcal L}^{E}(0,y)}{A_{\mathcal
L}^{E}(0,y)}.
\end{align*}

If $p>1$, the computation of the indices gives that
$$
I_{R_{l}^E}(\pi_E^* {\mathcal F},E_{red})  =  \frac{-H_{\mathcal
F}^{E}(c_{l}^{E})}{\prod_{j=1 \atop j \neq l}^{b_E}
(c_{l}^{E}-c_{j}^{E})}; \ \ I_{R_{l}^E}(\pi_E^* {\mathcal
L},E_{red}) =  \frac{-H_{\mathcal L}^{E}(c_{l}^{E})}{\prod_{j=1
\atop j \neq l}^{b_E} (c_{l}^{E}-c_{j}^{E})}
$$
and hence $H_{\mathcal F}^{E}(c_{l}^{E}) = H_{\mathcal
L}^{E}(c_{l}^{E})$ for $l=1,2,\ldots,b_E$. Consequently, we deduce
that the polynomials $H_{\mathcal F}^E(y)$ and $H_{\mathcal
L}^{E}(y)$ are equal.

Consider now the case $p=1$ which corresponds to $E=E_1$. We can
write
$$
H_{\mathcal F}^{E_1}(y)= a A_{\mathcal F}^\natural(y) + b
B_{\mathcal F}^{\natural}(y); \ \ \ H_{\mathcal L}^{E_1}(y)= a
A_{\mathcal L}^\natural(y) + b B_{\mathcal L}^{\natural}(y)
$$
with $A_{-}^{\natural}(y), B_{-}^{\natural}(y) \in {\mathbb C}[y]$.
Since $\pi_{E_1}$ is the blowing-up of the origin, it is easy to see
that
$$
A_{\mathcal F}^{E_1}(0,y)= A_{\mathcal F}^{r-1}(1,y) + y B_{\mathcal
F}^{r-1}(1,y); \ \ \  B_{\mathcal F}^{E_1}(0,y)=B_{\mathcal
F}^{r-1}(1,y)
$$
and similar equalities hold for the foliation $\mathcal L$. Thus,
from equation~\eqref{eq:A-E}, we deduce that
$$
A_{\mathcal F}^{\natural}(y) + y B_{\mathcal F}^{\natural}(0,y)=
A_{\mathcal L}^{\natural}(y) + y B_{\mathcal L}^{\natural}(0,y) =
\prod_{l=1}^{b_{E_1}}(y-c_{l}^{E_1}).
$$
Furthermore, the equality of the Camacho-Sad indices implies that
$B_{\mathcal F}^{\natural}(y) = B_{\mathcal L}^{\natural}(y)$ and
consequently $A_{\mathcal F}^{\natural}(y)=A_{\mathcal
L}^{\natural}(y)$. We conclude that $H_{\mathcal
F}^{E_1}(y)=H_{\mathcal L}^{E_1}(y)$ and this finish the proof of
the lemma.
\end{proof}
\begin{proof}[Proof of theorem~\ref{th:U_C}]
From the previous lemma we deduce that $\lambda \in U_{C}^{log}$ if
and only if, each foliation ${\mathcal F} \in {\mathbb
G}_{C,\lambda}^*$ is Zariski-general. This implies that
$U_{C}=U_{C}^{log}$ and the theorem follows straightforward.
\end{proof}
\begin{Remark}
Note that there are non Zariski-general foliations, even hamiltonian
ones. For instance, take $f=y(y-x^2)(2y-(1+\sqrt{-3})x^2)$ and
$\omega=df$; a generic polar curve of $\omega=0$ is irreducible with
one Puiseux pair equal to $(5,2)$ and hence the reduction of
singularities of $f=0$ is not a reduction of singularities of a
generic polar curve. Moreover, in this example $(1,1,1) \not \in
U_C$ whereas for $g=y(y-x^2)(y+x^2)$ a generic polar curve of $dg=0$
has two branches with coincidence equal to two and hence $(1,1,1)
\in U_C$. This shows  that the set $U_C$ depends on the analytic
type of the curve $C$.
\end{Remark}
\begin{Corollary}
If  ${\mathcal F} \in {\mathbb G}_{C,\lambda}^{\ast}$ is a
Zariski-general foliation, then the curves $C \cup \Gamma_{\mathcal
F}$ and $C \cup \Gamma_{{\mathcal L}_{\lambda}}$ are equisingular.
\end{Corollary}
Observe that the reciprocal of the corollary above is not true.
Consider  ${\mathcal F}$ defined by $\omega=0$ with $\omega = (4 i x
y^2 + 2 x^6y) dx + (y^2-2ix^2y-x^4-x^7)dy$. The foliation ${\mathcal
F}$ belongs to ${\mathbb G}_{C,\lambda}^{\ast}$ with
$C=(y(y-x^2)(y+x^2)=0)$ and $\lambda=(1,-i,i)$. The curves
$\Gamma_{\mathcal F}$ and $\Gamma_{{\mathcal L}_{\lambda}}$ are both
irreducible with one Puiseux pair equal to $(5,2)$. Hence $C \cup
\Gamma_{\mathcal F}$ and $C \cup \Gamma_{{\mathcal L}_{\lambda}}$
are equisingular. However, $\pi_C$ is not a reduction of
singularities of any of the generic polar curves and then $\lambda
\not \in U_C$. We also remark that ${\mathcal F}$ belongs to
${\mathbb G}_{C}^*$ although $(1,-i,i)$ is resonant.

\section{Kind equisingularity type}\label{sec:tipo-amable}
Let us consider a curve $C \subset ({\mathbb C}^2,0)$ which can have
singular branches and take $\rho:({\mathbb C}^2,0) \rightarrow
({\mathbb C}^2,0)$  any $C$-ramification. The existence of a curve
$Z$ such that $\rho^{-1}Z$ is a perfect adjoint curve of
$\rho^{-1}C$ can not be assured in general. We look for conditions
over $C$ that guarantee  the existence of perfect adjoint curves of
$\rho^{-1}C$ and, in this case, we  also define the equisingularity
type $\chi_C$.
\begin{Definition}
We say that a curve $C$ has a {\em kind equisingularity type\/} if
for each dead arc of $G(C)$ with bifurcation divisor $E_b$ and
terminal divisor $E_t$ we have that $m(E_b)=2 m(E_t)$.
\end{Definition}
Let us explain what having a kind equisingularity type means  in
terms of the equisingularity type of $C$. If $E_b$ is a bifurcation
divisor of $G(C)$ belonging to a dead arc with terminal divisor
$E_t$, then $m(E_b)=n_{E_b} m(E_t)$ by appendix~\ref{ap:grafodual}.
Hence, the curve $C$ has a kind equisingularity type if, and only
if, $n_{E_b}=2$ for each bifurcation divisor $E_b$ of $G(C)$ which
belongs to a dead arc. In particular, this implies that each dead
arc in $G(C)$ has only two vertices: the bifurcation divisor and the
terminal divisor. Observe that this property does not characterize
the fact of having a kind equisingularity type; it is enough to
consider the curve $y^3-x^5=0$ which does not have kind
equisingularity type. We have the following result of
characterization for  kind equisingularity types:
\begin{Proposition}\label{prop:exist-perfadj}
The following statements are equivalent:
\begin{itemize}
    \item The equisingularity type $\epsilon (C)$ is kind.
    \item There is a germ of curve $Z \subset ({\mathbb C}^2,0)$ such that
$\rho^{-1}Z$ is a perfect adjoint of $\rho^{-1}C$ for any
$C$-ramification $\rho$.
\end{itemize}
Moreover $\epsilon(C\cup Z)$ does not depend on  the choice of $Z$.
\end{Proposition}
\begin{proof}
Let $C \subset ({\mathbb C}^2,0)$ be a plane curve and consider
$\rho: ({\mathbb C}^2,0) \rightarrow ({\mathbb C}^2,0)$ any
$C$-ramification.

Assume first that there is a curve $Z$ such that $\rho^{-1}Z$ is a
perfect adjoint curve of $\rho^{-1}C$. Take any bifurcation divisor
$E$ of $G(C)$ which belongs to a dead arc with terminal divisor
$E_t$. Then $E$ is a Puiseux divisor and $m(E)=\underline{n}_E n_E$
with $n_E \geq 2$ and $m(E_t)=\underline{n}_E$. Let us  prove that
$n_E=2$.

Let $\{\tilde{E}^{j}\}_{j=1}^{\underline{n}_E}$ be the divisors
associated to $E$ in $G(\rho^{-1}C)$. We have that
\begin{equation}\label{eq:kind-1}
b_{\tilde{E}^{j}}= (b_E-1)n_E \ \ \mbox{ for all } \ \ j=1,\ldots,
\underline{n}_E.
\end{equation}
Let us denote by $b_{\tilde{E}^j}^*$ the number of  edges and arrows
which leave from $\tilde{E}^j$ in $G(\rho^{-1} C \cup \rho^{-1}Z)$.
Taking into account that $\rho^{-1}Z$ is a perfect adjoint of
$\rho^{-1}C$, from corollary~\ref{cor:perfadj-nosing} we have that
\begin{equation}\label{eq:kind-2}
b_{\tilde{E}^j}^* = 2 b_{\tilde{E}^j}-1 \ \ \mbox{ for all } \ \
j=1,\ldots, \underline{n}_E.
\end{equation}
 Moreover, using the
relationship between $G(C \cup Z)$ and $G(\rho^{-1}C \cup
\rho^{-1}Z)$, we can compute $b_{\tilde{E}^j}^*$ in terms of
$b_E^*$, where $b_E^*$ is the number of edges and arrows which leave
from $E$ in $G(C \cup Z)$. In fact, note that $E$ is also a Puiseux
divisor in $G(C \cup Z)$ and then there are two possibilities:
\begin{equation*}
b_{\tilde{E}^j}^* =
\left\{
  \begin{array}{ll}
    (b_E^*-1) n_E, & \hbox{if $E$ belong to a dead arc in $G(C\cup Z)$;} \\
    (b_E^*-1)n_E +1 , & \hbox{otherwise.}
  \end{array}
\right.
\end{equation*}
The first situation is not possible, because the equality
$b_{\tilde{E}^j}^* = (b_E^*-1) n_E$ and equations~\eqref{eq:kind-1},
\eqref{eq:kind-2} would imply that $2 n_E (b_E-1) -1 = (b_E^*-1)n_E$
and hence $n_E=1$ against the hypothesis. Then the second situation
holds so $b_{\tilde{E}^j}^* = (b_E^*-1) n_E+1$. Using again
equations~\eqref{eq:kind-1} and  \eqref{eq:kind-2}, we get that $(2
b_E - b^*_E -1) n_E=2$. Thus the only possible values are $n_E=2$
and $b^*_E = 2 b_E-2$.

Assume now that $C$ has a kind equisingularity type. Let $Z$ be a
plane curve such that $\pi_C$ gives a reduction of singularities of
$Z \cup C$ and that $G(C \cup Z)$ is obtained by adding to each
divisor $E$ of $G(C)$ the following number of arrows:
\begin{equation*}
\left\{
  \begin{array}{ll}
        b_E-1, & \hbox{if $E$ is a bifurcation divisor  which does not
        belong to a dead} \\
                & \hbox{arc in $G(C)$;} \\
       b_E-2, & \hbox{if $E$ is a bifurcation divisor  which belongs
to a dead arc in $G(C)$;} \\
       1, & \hbox{if $E$ is the terminal divisor of a dead arc in
$G(C)$;} \\
       0, & \hbox{in any other case.}
                  \end{array}
                \right.
\end{equation*}
Let us show that $\rho^{-1}Z$ is a perfect adjoint curve of
$\rho^{-1}C$. By the description of the reduction of singularities
of $Z$ given above, it is clear that $\rho^{-1}Z$ is composed only
by non-singular branches. We first prove that $\pi_{\rho^{-1}C}$
gives a reduction of singularities of $\rho^{-1} C \cup \rho^{-1}Z$.
Take any   branch $\gamma$ of $Z$ and consider the divisor $E$ of
$G(C)$ such that $\pi_{C}^*\gamma \cap E \neq \emptyset$. Let us see
that $\pi_{\rho^{-1}C}$ desingularizes $\rho^{-1}\gamma$. There are
three possible situations:
\begin{itemize}
    \item $E$ is a contact divisor with associated divisors
$\{\tilde{E}^j\}_{j=1}^{\underline{n}_E}$. Then $\rho^{-1} \gamma$
is composed by $\underline{n}_E$ non-singular branches and each of
them cuts one and only one divisor $\tilde{E}^j$.
    \item $E$ is a Puiseux divisor with associated divisors
$\{\tilde{E}^j\}_{j=1}^{\underline{n}_E}$. Then $\rho^{-1} \gamma$
is composed by $\underline{n}_E n_E$ non-singular branches and there
are exactly $n_E$ branches of $\rho^{-1} \gamma$ which cut each
$\tilde{E}^j$ in $n_E$ different points (see
appendix~\ref{ap:ramificacion}).
    \item $E$ is the extremity of a dead arc with bifurcation
divisor $E_b$. Let $\{\tilde{E}_b^j\}_{j=1}^{\underline{n}_{E_b}}$
be the divisors  associated to $E_b$. Then $\rho^{-1} \gamma$ is
composed by $\underline{n}_{E_b} = m(E)$ branches and each of them
cuts one and only one of the divisors $\tilde{E}_b^j$.
\end{itemize}
Moreover, $\pi_{\rho^{-1}C}$ is a reduction of singularities of
$\rho^{-1}Z$. In fact, consider two branches $\gamma$ and $\gamma'$
of $Z$ which cut the same divisor $E$ and let $\sigma$ and $\sigma'$
be two branches of $\rho^{-1} \gamma$ and $\rho^{-1}\gamma'$
respectively, such that they cut the same divisor $\tilde{E}^{j}$.
Then $\sigma$ and $\sigma'$ cut $\tilde{E}^j$ in different points
since otherwise the coincidence between $\gamma$ and $\gamma'$ would
be greater than $v(E)$. A similar argument proves that
$\pi_{\rho^{-1}C}$ is the minimal reduction of singularities of
$\rho^{-1}C \cup \rho^{-1}Z$.

In order to assure that $\rho^{-1}Z$ is a perfect adjoint of
$\rho^{-1}C$ we also need to check if $b_{\tilde{E}}^* = 2
b_{\tilde{E}}-1$ for each bifurcation divisor $\tilde{E}$ of
$G(\rho^{-1}C)$. Let $E$ be the bifurcation divisor of $G(C)$ which
$\tilde{E}$ is associated to. Let us consider the three possible
cases for $E$:
\begin{itemize}
    \item $E$ is a contact divisor in $G(C)$ and we have that
$b_{\tilde{E}}=b_E$ and $b_{E}^* = 2 b_E-1$. But $E$ is also a
contact divisor in $G(C \cup Z)$ and hence $b_{\tilde{E}}^* =
b_{E}^*$. We deduce that $b_{\tilde{E}}^*=2 b_{\tilde{E}}-1$.
    \item $E$ is a Puiseux divisor belonging to a dead arc in $G(C)$ and hence
$b_{\tilde{E}}= (b_E-1)n_E$ and $b_{E}^*= 2b_E-2$. In this case, $E$
is a Puiseux divisor without dead arc in $G(C \cup Z)$ and we have
that $b_{\tilde{E}}^* = (b_E^*-1)n_E+1$. We deduce that
$b_{\tilde{E}}^* = 2 b_{\tilde{E}}-n_E+1$ and the result follows
since by hypothesis $n_E=2$.
    \item $E$ is a Puiseux divisor without a dead arc in $G(C)$,
thus $b_{\tilde{E}}=(b_E-1)n_E +1$ and $b_{E}^* = 2 b_E -1$. The
divisor $E$ is also a Puiseux divisor without a dead arc in $G(C
\cup Z)$, so $b_{\tilde{E}}^*=(b_E^*-1)n_E +1$. Hence we conclude
that $b_{\tilde{E}}^*= 2 b_{\tilde{E}}-1$.
\end{itemize}
It is clear that the equisingularity type $\epsilon(C \cup Z)$ does
not depend on the choice of the curve $Z$.
\end{proof}
If $C$ is a curve with kind equisingularity type, we say that $Z$ is
a {\em perfect adjoint curve\/} of $C$ if $\rho^{-1}Z$ is a perfect
adjoint curve of $\rho^{-1}C$, for any $C$-ramification $\rho$. We
are interested in the description of the equisingularity type
$\chi_C = \epsilon(C \cup Z)$. A first result in this direction is
the following lemma:
\begin{Lemma}\label{lem:coinc-perfadj}
Consider  a curve  $C$  with kind equisingularity type and let $Z$
be a perfect adjoint curve of $C$ with  $Z= \cup_{E \in B(C)}Z^E$.
Then ${\mathcal C}(\zeta^E,\xi^E)=v(E)$ for any two branches
$\zeta^E, \xi^E$ of $Z^E$.
\begin{proof}
The result follows from corollary~\ref{cor:perfadj-nosing} and
equation~\eqref{eq:coincidencias-ram}.
\end{proof}
\end{Lemma}
The next proposition gives a completely description of
$\chi_C=\epsilon(C\cup Z)$ in terms of $\epsilon(C)$:
\begin{Proposition}\label{prop:equis-adjperf}
Let $C$ be a curve with kind equisingularity type and $Z$ a perfect
adjoint curve of $C$. Then $\pi_C$ gives a reduction of
singularities of $Z\cup C$. Moreover, the branches of $Z$ intersect
an irreducible component $E$ of the exceptional divisor of $\pi_C$
as follows:
\begin{itemize}
    \item If $E$ is a bifurcation divisor of $G(C)$,  the number
    of branches of $Z$ cutting $E$
    equals to $b_E-2$ if $E$ is in a dead arc and to
    $b_E-1$ otherwise.
    \item If $E$ is a terminal divisor of a dead arc of $G(C)$, there is
    exactly one branch of $Z$
    through $E$.
\item Otherwise, no branches of $Z$  intersect $E$.
\end{itemize}
\end{Proposition}
Remark that the fact that ``$\pi_C$ gives a reduction of
singularities of $C \cup Z$" does not imply that $\pi_{\rho^{-1}C}$
desingularizes $\rho^{-1} C \cup \rho^{-1} Z$. However, the
description of the dual graph $G(C \cup Z)$ given in
proposition~\ref{prop:equis-adjperf} characterizes the fact of $Z$
being a perfect adjoint curve of $C$ whenever $C$ has a kind
equisingularity type. In fact, in
proposition~\ref{prop:exist-perfadj} we have already proved that, if
$C$ has a kind equisingularity type,  a curve $Z$ such that $G(C
\cup Z)$ is as described in proposition~\ref{prop:equis-adjperf} is
a perfect adjoint curve of $C$ and the proof of
proposition~\ref{prop:equis-adjperf}  will show the reciprocal.

In order to prove proposition~\ref{prop:equis-adjperf} we first
describe the equisingularity type of the irreducible components of
$Z$ in terms of the equisingularity data of $C = \cup_{i=1}^{r}
C_i$. Given an irreducible component $C_i$ of $C$ we denote by
$\{\beta_{0}^{i},\beta_{1}^{i}, \ldots, \beta_{g_i}^{i}\}$ its
charac\-teristic exponents, $\{(m_j^{i},n_{j}^{i})\}_{j=1}^{g_i}$
the Puiseux pairs of $C_i$ and $n^i$ is the multiplicity $m_0(C_i)$
at the origin. We use  the notations introduced in
appendix~\ref{ap:grafodual} for the dual graph $G(C)$.

\begin{Lemma}\label{lema:equising-perfadj}
Consider a curve $C$ with kind equisingularity type and let $Z$ be
perfect adjoint curve of $C$ with decomposition $Z=\cup_{E \in B(C)}
Z^{E}$. Then, for each $E \in B(C)$ and $i \in I_{E}^*$, we have
that
\begin{itemize}
    \item[(i)]  If $E$ is a contact divisor, then the curve
$Z^{E}$ has $b_E-1$ irreducible components. Each irreducible
component $\zeta$ of $Z$ with characteristic exponents $\{
\nu_0^{\zeta}, \nu_{1}^{\zeta}, \ldots, \nu_{k_E}^{\zeta}\}$ given
by
$$
\nu_0^{\zeta} = m_0(\zeta)=\underline{n}_E, \ \ \nu_{l}^{\zeta}=
\underline{n}_E \beta_{l}^{i}/n^{i} \text{ for } l=1,2, \ldots, k_E.
$$
    \item[(ii)] If $E$ is a Puiseux divisor which belongs to a dead
arc, the curve $Z^E$ has one irreducible component $\zeta_0$ with
characteristic exponents \linebreak $\{ \nu_0^{\zeta_0},
\nu_{1}^{\zeta_0}, \ldots, \nu_{k_E}^{\zeta_0}\}$ given by
$$
\nu_0^{\zeta_0} = m_0(\zeta_0)=\underline{n}_E, \ \
\nu_{l}^{\zeta_0}= \underline{n}_E \beta_{l}^{i}/n^{i} \text{ for }
l=1,2, \ldots, k_E,
$$
and $b_{E}-2$ irreducible components such that each branch
$\zeta
\subset Z^E\smallsetminus \zeta_0$ has characteristic exponents $\{
\nu_0^{\zeta}, \nu_{1}^{\zeta}, \ldots, \nu_{k_E}^\zeta,
\nu_{k_E+1}^{\zeta}\}$ given by
$$
\nu_0^{\zeta} = m_0(\zeta)=\underline{n}_E n_E, \ \ \nu_{l}^{\zeta}=
\underline{n}_E n_E \beta_{l}^{i}/n^{i} \text{ for } l=1,2, \ldots,
k_E+1.
$$
    \item[(iii)] If $E$ is a bifurcation divisor which does not
belong to a dead arc, then $Z^E$ has $b_E-1$ irreducible components.
Each irreducible component $\zeta$ of $Z$ with characteristic
exponents $\{ \nu_0^{\zeta}, \nu_{1}^{\zeta},
\ldots,\nu_{k_E}^\zeta, \nu_{k_E+1}^{\zeta}\}$ given by
$$
\nu_0^{\zeta} = m_0(\zeta)=\underline{n}_E n_E, \ \ \nu_{l}^{\zeta}=
\underline{n}_E n_E\beta_{l}^{i}/n^{i} \text{ for } l=1,2, \ldots,
k_E+1.
$$
\end{itemize}
\end{Lemma}
\begin{proof}
Consider any $C$-ramification $\rho:({\mathbb C}^2,0) \rightarrow
({\mathbb C}^2,0)$ and denote $\tilde{C}=\rho^{-1}C$. Let
$\{\tilde{E}^{l} \}_{l=1}^{\underline{n}_E}$ be the divisors of
$G(\tilde{C})$ associated to a divisor $E$ of $G(C)$. By the results
in section~\ref{sec:adjuntos}, we have that $\rho^{-1}Z^E =
\cup_{j=1}^{\underline{n}_E} \tilde{Z}^{\tilde{E}^j}$ where
$\tilde{Z}=\cup_{\tilde{E} \in B(\tilde{C})} \tilde{Z}^{\tilde{E}}$
is the decomposition of $\tilde{Z}=\rho^{-1}Z$. Let us study the
different possibilities for $E$:

\vspace{\baselineskip} \hspace{.5cm}(i) {\em E is a contact
divisor:\/} then $v(E)=m_E /\underline{n}_E$ with $m_{E} >
m_{k_E}^{i}$ and $\underline{n}_E=n_{1}^{i} \cdots n_{k_E}^{i}$ for
any $i \in I_E$. Consequently, the $k_E$ first Puiseux pairs of an
irreducible component $\zeta^E$ of $Z^E$ coincide with the ones of
$C_i$, for any $i \in I_E$, since ${\mathcal
C}(\zeta^{E},C_i)=v(E)$. Thus, a Puiseux series of $\zeta^E$ is
given by
$$
\varphi_{\zeta}(x) = \sum_{i < \tau} a_i x^i + a_{\tau} x^{\tau} +
\cdots + a^{\zeta} x^{v(E)} + \cdots,
$$
where $\tau=m_{k_E}^{i} /\underline{n}_E$ and $a_{\tau} \neq 0$.
This implies that $m_0(\zeta^E)=d \cdot \underline{n}_E$. Let us
show that $m_{0} (\zeta^E) = \underline{n}_E$.

We have that $\tilde{\zeta}^E=\rho^{-1} \zeta^E \subset
\cup_{l=1}^{\underline{n}_E} \tilde{Z}^{\tilde{E}_l}$, and if we
write $\tilde{\zeta}^E =
\cup_{l=1}^{\underline{n}_E}\tilde{\zeta}^{\tilde{E}_l}$ with
$\tilde{\zeta}^{\tilde{E}_l} \subset \tilde{Z}^{\tilde{E}_l}$, then
$m_0(\tilde{\zeta}^{\tilde{E}_l}) \geq 1$. By
corollary~\ref{cor:perfadj-nosing}, each curve
$\tilde{\zeta}^{\tilde{E}_l}$ has $m_0(\tilde{\zeta}^{\tilde{E}_l})$
non-singular irreducible components and the coincidence between two
of them is equal to $v(\tilde{E}_l)$. Moreover, the irreducible
components of $\tilde{\zeta}^{E}$ are in bijective correspondence
with the Puiseux series of $\zeta^E$. Then, if $a^\zeta \neq 0$, the
coefficients of $x^{v(E)}$ in the different Puiseux series of
$\zeta^{E}$ are given by $a^{\zeta} \xi^{v(E) m_0(\zeta^E)}$ with
$\xi^{m_0(\zeta^E)}=1$. But since
$$
v(E) \cdot m_0(\zeta^E) = \frac{m_E}{\underline{n}_E} \cdot
m_0(\zeta^E) = m_E \cdot d
$$
then $a^{\zeta} \xi^{v(E) m_0(\zeta^E)}$ takes at most
$\underline{n}_E$ different values and hence $d=1$. If $a^\zeta =0$,
then $m_0(\zeta^E)=\underline{n}_E$ since otherwise one of the
curves $\tilde{\zeta}^{\tilde{E}_l}$ has at least two irreducible
components with coincidence greater than $v(\tilde{E}_l)$.

We deduce that each irreducible component $\zeta^E$ of $Z^E$ has
multiplicity equal to $\underline{n}_E$. Since $m_0(Z^E)=
\underline{n}_E (b_{E}-1)$, then $Z^E$ has exactly $b_E-1$
irreducible components with multiplicity $\underline{n}_E$.
Moreover, the Puiseux pairs of each irreducible component $\zeta$ of
$Z^E$ coincide with the $k_E$ first Puiseux pairs of $C_i$ for $i
\in I_E$ and the characteristic exponents $\{ \nu_0^{\zeta},
\nu_{1}^{\zeta}, \ldots, \nu_{k_E}^{\zeta}\}$ of $\zeta$ are given
by $\nu_{l}^{\zeta} = \underline{n}_E \beta_{l}^{i} /n^i$ for
$l=0,1,\ldots, k_E$.

\vspace{\baselineskip} \hspace{.5cm}(ii) {\em $E$ is a Puiseux
divisor which belongs to a dead arc:\/} we have that $v(E)=m_E
/\underline{n}_E n_E$ with $n_E =2$ because $C$ has a kind
equisingularity type and then $m_0(Z^{E})=\underline{n}_E (n_E
(b_E-1)-1) = \underline{n}_E n_E (b_E-2) + \underline{n}_E$.

An irreducible component $\zeta^E$ of $Z^E$ has at least the $k_E$
first Puiseux pairs equal to the ones of $C_i$ with $i \in I_E$.
Thus $m_0(\zeta^E) \geq \underline{n}_E$. A Puiseux series
$\varphi_\zeta (x)$ of $\zeta^E$ is given by
$$
\varphi_{\zeta}(x)= \sum_{l < v(E)} a_l x^l + a^\zeta x^{v(E)} +
\ldots,
$$
but since $\underline{n}_En_E$ does not divide $m_0(Z^E)$, then
there is at least one irreducible component $\zeta_0^E$ of $Z^E$
such that the coefficient $a^{\zeta_0}$ of $x^{v(E)}$ is zero.
Moreover, $\zeta_0^E$ must be unique because the existence of
another irreducible component $\delta_0^E$  of $Z^E$ with
$a^{\delta_0} =0$ would imply that  ${\mathcal C}(\zeta_0^E,
\delta_0^E)
>v(E)$ in contradiction with lemma~\ref{lem:coinc-perfadj}. Let us show
that $m_0(\zeta_0^E)=\underline{n}_E$. In fact, $m_0(\zeta_0^E)= d
\cdot \underline{n}_E$ with $d \in {\mathbb N}$. Consider the curve
$\tilde{\zeta}_0^E=\rho^{-1} \zeta_0^E$ and write
$\tilde{\zeta}_0^{E}= \cup_{l=1}^{\underline{n}_E}
\tilde{\zeta}^{\tilde{E}_l}_0$ with $\tilde{\zeta}^{\tilde{E}_l}_0
\subset \tilde{Z}^{\tilde{E}_l}$. By
corollary~\ref{cor:perfadj-nosing}, the number of irreducible
components of $\tilde{Z}^{\tilde{E}_l}$ is equal to its
multiplicity, hence $m_0(\tilde{\zeta}^{\tilde{E}_l}_0)=1$ since
otherwise the coincidence between two branches of
$\tilde{\zeta}^{\tilde{E}_l}_0$ will be greater than
$v(\tilde{E}_l)$. Hence $m_0(\zeta^E_0)=\underline{n}_E$.
Consequently, we have that
$$
m_0(Z^E \smallsetminus \zeta_0^E) = \underline{n}_E n_E (b_E-2).
$$
Consider now  an irreducible component $\zeta^E$  of $Z^E
\smallsetminus \zeta^E_0$. The coefficient $a^\zeta$ in
$\varphi_\zeta(x)$ must be non-zero and thus $m_0(\zeta^E) \geq
\underline{n}_E n_E$. With similar arguments as above, we show that
$m_0(\zeta^E)= \underline{n}_E n_E$.

We have proved that $Z^E$ has one irreducible component $\zeta_0^E$
with multiplicity $\underline{n}_E$ and $b_E-2$ irreducible
components  with multiplicity $\underline{n}_E n_E$. The
characteristic exponents $\{ \nu_0^{\zeta_0}, \nu_{1}^{\zeta_0},
\ldots, \nu_{k_E}^{\zeta_0}\}$ of $\zeta_0^E$ are given by
$\nu_{l}^{\zeta_0} = \underline{n}_E \beta_{l}^{i} /n^i$, for
$l=1,\ldots,k_E$, and the characteristic exponents $\{
\nu_0^{\zeta}, \nu_{1}^{\zeta}, \ldots, \nu_{k_E+1}^{\zeta}\}$  of a
branch $\zeta^E$ of  $Z^E \smallsetminus \zeta^E_0$ are given by
$\nu_{l}^{\zeta} = \underline{n}_E n_E \beta_{l}^{i} /n^i$ for
$l=0,1,\ldots, k_E+1$ and $i \in I_E$.

 \vspace{\baselineskip} \hspace{.5cm}(iii) {\em $E$ is a Puiseux
divisor which does not belong to a dead arc:\/} we have that
$v(E)=m_E/\underline{n}_E n_E$ with $n_E >1$. Take any irreducible
component $\zeta^E$ of $Z^E$. Let us see that $m_0(\zeta^E)=
\underline{n}_E n_E$. Consider
$$
\varphi_{\zeta}(x) = \sum_{l < v(E)} a_l x^l + a^{\zeta} x^{v(E)} +
\cdots
$$
a Puiseux series of $\zeta^E$. The hypothesis over $E$ imply that
$(m_E,n_E)$ is not a Puiseux pair of $C_j$ if $j \in I_E
\smallsetminus I_E^*$, or equivalently, the coefficient of
$x^{v(E)}$ in the Puiseux series of $C_j$ is zero.  In particular,
we deduce that $a^\zeta \neq 0$ for all irreducible components
$\zeta^E$ of $Z^E$ since ${\mathcal C}(C_j,\zeta^E)=v(E)$.
Consequently, $(m_E,n_E)$ is a Puiseux pair of $\zeta^E$ and the
$k_E+1$ Puiseux pairs of $\zeta^E$ coincide with the ones of $C_i$
with $i \in I_E^*$. With similar arguments as in case (i) we prove
that $m_0(\zeta^E)=\underline{n}_E n_E$.

From the fact that $m_0(Z^E)=\underline{n}_E n_E (b_E-1)$, we deduce
that $Z^E$ has exactly $b_E-1$ irreducible components, each of them
with multiplicity $\underline{n}_E n_E$. Hence, the characteristic
exponents $\{ \nu_0^{\zeta}, \nu_{1}^{\zeta}, \ldots,
\nu_{k_E+1}^{\zeta}\}$ of a branch $\zeta^E$ of $Z^E$ are given by
$\nu_{l}^{\zeta} = \underline{n}_E n_E \beta_{l}^{i}/n^i$ for
$l=1,\ldots, k_E+1$ and $i \in I_E^*$.
\end{proof}
The previous description  of the equisingularity type of the
irreducible components of $Z^E$ will be useful in the proof of
proposition~\ref{prop:equis-adjperf}.
\begin{proof}[Proof of proposition~\ref{prop:equis-adjperf}]
Let $C$ be a curve with kind equisingularity type and let $\pi_C : M
\rightarrow ({\mathbb C}^2,0)$ be its minimal reduction of
singularities. Consider $Z$ a perfect adjoint curve of $Z$  with
decomposition $Z = \cup_{E \in B(C)} Z^E$ satisfying properties
D1.-D5. in section~\ref{sec:adjuntos}. It is clear that the points
of $\pi_C^* Z \cap \pi^{-1}_C(0)$ coincide with the union of the
sets $\pi_{C}^{\ast} Z^E \cap \pi_{C}^{-1}(0)$ for $E \in B(C)$. We
deduce that if $Z$ cuts a divisor $E$, then $E$ is either a
bifurcation divisor or it belongs to a dead arc, but since each dead
arc of $G(C)$ has only to vertices, then $E$ is either a bifurcation
or a terminal divisor.

Assume first that $E$ is a bifurcation divisor without a dead arc
attached to it. Then properties D3.-D5. of the decomposition of $Z$
imply that each irreducible component $\zeta^E$ of $Z^E$ cuts $E$,
i.e., $\pi_{E}^* \zeta^E \cap E_{red} \neq \emptyset$. Moreover, the
number of points of $\pi_{E}^* Z^E \cap E_{red}$ is equal to the
number of irreducible components of $Z^E$. In fact,  if $\pi_{E}^*
\zeta^E \cap E_{red} = \pi_{E}^* \xi^E \cap E_{red}$ then ${\mathcal
C}(\zeta^E,\xi^E) >v(E)$ in contradiction with
lemma~\ref{lem:coinc-perfadj}. The present hypothesis correspond to
the cases (i) and (iii) of lemma~\ref{lema:equising-perfadj}, hence
the number of points of $\pi_{E}^* Z^E \cap E_{red}$ is equal to
$b_E-1$. It is clear that $\pi_E$ is a reduction of singularities of
each irreducible component $\zeta^E$ of $Z^E$ since each curve
$\pi_{E}^*\zeta^E$ is an  $E_{red}$-curvette by
lemma~\ref{lema:equising-perfadj}.

Assume now that $E$ is a bifurcation divisor which belong to  a dead
with terminal divisor $E_t$. By properties D3.-D5. of the
decomposition of $Z$, we have that either $\pi_{E}^* \zeta^E \cap
E_{red} \neq \emptyset$ or $\pi_{E}^* \zeta^E \cap \pi_{E}'(E_t)
\neq \emptyset$ for an irreducible component $\zeta^E$ of $Z^E$. By
lemma~\ref{lema:equising-perfadj}, there is an irreducible component
$\zeta_0^E$ of $Z^E$ with multiplicity $\underline{n}_E$, thus
$\pi_E^* \zeta_{0}^E \cap \pi_{E}'(E_t) \neq \emptyset$ since each
curve $\gamma$ with $\pi_E^*\gamma \cap E_{red} \neq \emptyset$ must
have multiplicity $\geq m(E)=\underline{n}_E n_E$. Moreover,
$\zeta_0^E$ is the only irreducible component of $Z^E$ which cuts
$E_t$ because the existence of another one $\xi^E_0$ would imply
that ${\mathcal C}(\zeta_0^E,\xi^E_0) \geq v(E_t) >v(E)$ in
contradiction with lemma~\ref{lem:coinc-perfadj}. Finally, the
number of points of $\pi_{E}^* Z^E \cap E_{red}$ coincides with the
number of irreducible components of $Z^E \smallsetminus \zeta_0^E$
which is $b_E -2$. We also have that $\pi_E$ is a reduction of
singularities of $Z^E$ since $\zeta_0^E$ is a $\pi_{E}'
(E_t)$-curvette
 and $\zeta^E$ is an $E_{red}$-curvette
for each $\zeta^E \subset Z^E \smallsetminus \zeta_0^E$ by
lemma~\ref{lema:equising-perfadj}.

The fact that $\pi_C$ gives a reduction of singularities of $C \cup
Z$ follows using property D2. and the result is proved.
\end{proof}

\section{Proof of the main theorem}\label{sec:main-th}

Consider a curve $C=\cup_{i=1}^{r}C_i$ which can have singular
branches. Let $U_C$ be the set of $\lambda \in {\mathbb P}_{\mathbb
C}^{r-1}$ such that there exists ${\mathcal F} \in {\mathbb
G}_{C,\lambda}^*$ with $\rho^{-1}\Gamma_{\mathcal F}$ a perfect
adjoint curve of $\rho^{-1}C$, for any $C$-ramification $\rho$. This
section is devoted to prove the following result:
\begin{Theorem}\label{th:U-C-no-vacio}
The  set $U_C$ is a non-empty Zariski open set if and only if $C$
has a kind equisingularity type. Moreover, in this case
$\wp({\mathcal F})=\chi_C$ for any ${\mathcal F} \in {\mathbb
G}_C^*$ with $\lambda({\mathcal F}) \in U_C$.
\end{Theorem}

Take any $C$-ramification $\rho: ({\mathbb C}^2,0) \rightarrow
({\mathbb C}^2,0)$ given by $x=u^n, y=v$. Consider a foliation
${\mathcal F} \in {\mathbb G}_{C,\lambda}^{\ast}$, then the
transform $\rho^*{\mathcal F}$ belongs to
$G_{\rho^{-1}C,\lambda^*}^{\ast}$ where $\lambda^* =
\lambda(\rho^*{\mathcal F}) \in {\mathbb P}_{\mathbb C}^{m-1}$ and
$m=m_0(C)$ is the multiplicity of $C$ at the origin. We denote by
$\Gamma_{\mathcal F}$ and $\Gamma_{\rho^{\ast} {\mathcal F}}$ two
generic polar curves of $\mathcal F$ and $\rho^* {\mathcal F}$
respectively.

It is clear that the foliation $\rho^{\ast}{\mathcal F}$ has a curve
of separatrices with only non-singular branches. Consequently, by
the results of section~\ref{sec:non-sp}, $\Gamma_{\rho^*{\mathcal
F}}$ is a perfect adjoint curve of $\rho^{-1}C$ if and only if
$\lambda^* \in U_{\rho^{-1}C}$ and in that case,
$\epsilon(\Gamma_{\rho^* {\mathcal F}} \cup
\rho^{-1}C)=\chi_{\rho^{-1}C}$. However, in general,
$\rho^{-1}\Gamma_{\mathcal F}$ and $\Gamma_{\rho^{\ast} {\mathcal
F}}$ are not equi\-singular (see~\cite{Cor-03}). Consider the
following properties:
\begin{align*}
(A): & \ \  \epsilon(\Gamma_{\rho^* {\mathcal F}} \cap \rho^{-1} C)
= \chi_{\rho^{-1}C} \\
(B) : & \ \ \epsilon(\rho^{-1} \Gamma_{\mathcal F} \cap \rho^{-1}C)
= \chi_{\rho^{-1}C}
\end{align*}
\begin{Proposition}\label{prop:A-B}
Property $(A)$ implies $(B)$. Moreover, both properties are
equivalent if the curve $C$ has at most two different tangent lines.
\end{Proposition}
Observe that properties $(A)$ and $(B)$ above do not depend on the
choice of the $C$-ramification $\rho$.
\begin{Definition}
We say that $\mathcal F$ is a {\em Zariski-general foliation\/} when
property $(B)$ holds.
\end{Definition}

\begin{Notation}
In this section, we denote by $\tilde{C}$ and $\tilde{\Gamma}$ the
curves $\rho^{-1}C$ and  $\rho^{-1} \Gamma_{\mathcal F}$
respectively; the transform of the polar $\rho^{-1} \Gamma({\mathcal
F};[a:b])$ will be denoted by $\tilde{\Gamma}_{[a:b]}$ or
$\tilde{\Gamma}_{\mathcal F}$ when the explicit direction of
polarity or the foliation are needed. If $\pi_{\tilde{C}}: \tilde{M}
\rightarrow ({\mathbb C}^2,0)$ is the minimal reduction of
singularities of $\tilde{C}$, we denote by $\tilde{E}$ an
irreducible component of $\pi_{\tilde{C}}^{-1}(0)$ and by
$\pi_{\tilde{E}}: \tilde{M}_{\tilde{E}} \rightarrow ({\mathbb
C}^2,0)$ the morphism reduction of $\pi_{\tilde{C}}$ to $\tilde{E}$.
The reader could refer to appendix~\ref{ap:ramificacion} for a
detailed description of the ramification tools.
\end{Notation}
Let us state two lemmas concerning the infinitely near points of
$\tilde{\Gamma}$ and $\Gamma_{\rho^* {\mathcal F}}$.
\begin{Lemma}\label{lema-ptosE1}
Consider a foliation ${\mathcal F} \in {\mathbb G}_{C}^*$ and let
$\tilde{E}_1$ be the irreducible component  of
$\pi_{\tilde{C}}^{-1}(0)$ with $v(\tilde{E}_1)=n$. Then the set
$$
\pi_{\tilde{E}_1}^* \tilde{\Gamma}_{[a:b]} \cap \tilde{E}_{1,red}
\smallsetminus \pi_{\tilde{E}_1}^* \tilde{C} \cap \tilde{E}_{1,red}
$$
has exactly $b_{\tilde{E}_{1}}-1$ points which depend on $[a:b]$.
\end{Lemma}
\begin{proof}
Observe that the divisor $\tilde{E}_1$ of $\pi_{\tilde{C}}^{-1}(0)$
is associated to the divisor $E_1$ of $\pi_{C}^{-1}(0)$ and hence
the coordinates $(x,y)$ and $(u,v)$ are adapted to $E_1$ and
$\tilde{E}_1$, respectively. Let $\omega=A(x,y) dx + B(x,y)dy$ be a
1-form defining $\mathcal F$. Then $\Gamma_{[a:b]}$ is defined by $a
A(x,y) + b B(x,y)=0$ and $\tilde{\Gamma}_{[a:b]}$ is given by $a
A(u^n,v) + b B(u^n,v)=0$. Take coordinates $(\tilde{u},\tilde{v})$
in the first chart of $\tilde{E}_1$ such that
$\pi_{\tilde{E}_1}(\tilde{u},\tilde{v})=(\tilde{u},\tilde{u}^n
\tilde{v})$ and $\tilde{E}_1=(\tilde{u}=0)$. The strict transform
$\pi_{\tilde{E}_1}^* \tilde{\Gamma}_{[a:b]}$ is given by
$$
\pi_{\tilde{E}_1}^* \tilde{\Gamma}_{[a:b]}= \{ a A_{\nu}
(1,\tilde{v}) + b B_{\nu}(1,\tilde{v}) + \tilde{u} ( \cdots ) =0 \},
$$
where $\nu=\nu_0({\mathcal F})$ and $A_{\nu}(x,y) dx +
B_{\nu}(x,y)dy$ is the $\nu$-jet of $\omega$. Then the points of
$\pi_{\tilde{E}_1}^* \tilde{\Gamma}_{[a:b]} \cap \tilde{E}_{1,red}$
are defined by $\tilde{u}=0$ and $a A_{\nu}(1,\tilde{v}) + b
B_{\nu}(1,\tilde{v}) =0$. Taking into account that
$\tilde{\Gamma}_{[a:b]}$ is a strict adjoint of $\tilde{C}$ and
using similar arguments as in the proof of
proposition~\ref{prop:U_C-log} case $p=1$, we get that the points of
$\pi_{\tilde{E}_1}^* \tilde{\Gamma}_{[a:b]} \cap \tilde{E}_{1,red}
\smallsetminus \pi_{\tilde{E}_1}^* \tilde{C} \cap \tilde{E}_{1,red}$
are given by $\tilde{u}=0$ and $H^{\tilde{E}_1}(\tilde{v})=0$ with
$$
H^{\tilde{E}_1} (v) = a A_{\nu}^{\natural}(v) + b B_{\nu}^{\natural}
(v),
$$
where $A_{\nu}^{\natural}(v)$ and $B_{\nu}^{\natural}(v)$ do not
have common roots. Thus the result follows straightforward.
\end{proof}
\begin{Corollary}
Given a foliation ${\mathcal F} \in {\mathbb G}_{C}^*$, the set
$\pi_{E_1}^* \Gamma_{[a:b]}^{\mathcal F} \cap E_{1,red}
\smallsetminus \pi_{E_1}^{\ast}C \cap E_{1,red}$ has exactly
$b_{E_{1}}-1$ points which depend on $[a:b]$.
\end{Corollary}
\begin{proof}
The result follows from the fact that there is a bijection between
the points in $E_{1,red}$ and the ones in $\tilde{E}_{1,red}$ (see
lemma~\ref{lema:apB}).
\end{proof}

\begin{Lemma}\label{lem-polar-ram}
Consider a foliation ${\mathcal F} \in {\mathbb G}_{C}^*$. Then we
have that
$$
\pi_{\tilde{E}}^{\ast} \tilde{\Gamma} \cap \tilde{E}_{red} =
\pi_{\tilde{E}}^{\ast} \Gamma_{\rho^{\ast}{\mathcal F}} \cap
\tilde{E}_{red}
$$
for each irreducible component $\tilde{E}$ of
$\pi_{\tilde{C}}^{-1}(0)$ with $v(\tilde{E}) >n$. Moreover, $
m_{P}(\pi_{\tilde{E}}^{\ast} \tilde{\Gamma}) =
m_P(\pi_{\tilde{E}}^{\ast} \Gamma_{\rho^{\ast}{\mathcal F}})$ for
each $P \in \pi_{\tilde{E}}^{\ast} \tilde{\Gamma} \cap
\tilde{E}_{red}$.
\end{Lemma}

\begin{proof}
Let $\omega= A(x,y) dx + B(x,y) dy$ be a 1-form defining $\mathcal
F$. Then the curves $\tilde{\Gamma}$ and $\Gamma_{\rho^* {\mathcal
F}}$ are given by
\begin{align*}
\tilde{\Gamma} & = \{ a A(u^n,v) + b B(u^n,v)=0 \}; \\
\Gamma_{\rho^* {\mathcal F}} & = \{ a A(u^n,v) n u^{n-1} + b
B(u^n,v) =0 \}
\end{align*}
Take any irreducible component $\tilde{E}$ of
$\pi_{\tilde{C}}^{-1}(0)$ with $v(\tilde{E})=p >n$ and assume that
$(u,v)$ are coordinates adapted to $\tilde{E}$. By the results in
section~\ref{sec:loc-inv}, it is enough to prove that
\begin{equation}\label{eq:In-A-B}
In_p(a A^* + b \tilde{B};u,v) = In_p(a \tilde{A}+b\tilde{B};u,v) =
In_p (b \tilde{B};u,v)
\end{equation}
where $\tilde{A}(u,v)=n u^{n-1} A(u^n,v)$, $\tilde{B}(u,v)=B(u^n,v)$
and $A^* (u,v) = A(u^n,v)$.

Let $i + pj =k$ be the equation of the line which contains the side
of ${\mathcal N}(\rho^* {\mathcal F};u,v)$ with slope equal to
$-1/p$. Then it is clear that $\Delta(\rho^* \omega) \subset \{
(i,j) \in {\mathbb R}^2 \ : \ i +pj \geq k \}$. Moreover, $\Delta(a
\tilde{A} + b \tilde{B}) \subset \{ (i,j) \ : \ i+pj \geq k-p\}$ by
lemma~\ref{lem:polNew:F-Gamma}. Let us prove that
$\Delta(\tilde{A})$ and $\Delta(A^*)$ are contained in $\{(i,j) \ :
\ i + pj  > k-p\}$. Consider two cases:
\begin{itemize}
    \item If $(i,j)  \in \Delta(\tilde{A})$ then $(i+1,j) \in
    \Delta(\rho^* \omega)$ and hence $i + pj \geq k-1 > k-p$.
    \item If $(i,j) \in \Delta(A^*)$ then $(i+n,j) \in \Delta(\rho^*
    \omega)$ and consequently $i + pj \geq k-n > k-p$.
\end{itemize}
Thus the equalities in \eqref{eq:In-A-B} hold and the lemma is
proved.
\end{proof}

Let us show now that being a Zariski-general foliation only depends
on $\lambda({\mathcal F})$.
\begin{Proposition}
A foliation ${\mathcal F} \in {\mathbb G}_{C,\lambda}^{\ast}$ is
Zariski-general  if and only if ${\mathcal L}_{\lambda}$ is a
Zariski-general foliation.
\end{Proposition}
\begin{proof}
Let $\Gamma_{\mathcal F}$ and $\Gamma_{\mathcal L}$ be  generic
polar curves of $\mathcal F$ and ${\mathcal L}={\mathcal
L}_{\lambda}$, respectively, and put $\tilde{\Gamma}_{\mathcal
F}=\rho^{-1}\Gamma_{\mathcal F}$ and $\tilde{\Gamma}_{\mathcal
L}=\rho^{-1}\Gamma_{\mathcal L}$. Let us prove that  the infinitely
near points of $\tilde{\Gamma}_{\mathcal F}$ and
$\tilde{\Gamma}_{\mathcal L}$ coincide at each irreducible component
$\tilde{E}$ of $\pi_{\tilde{C}}^{-1}(0)$, $\tilde{E} \neq
\tilde{E}_1$. In fact, by lemma~\ref{lema-puntos-polarF-L}, we have
that
$$
\pi_{\tilde{E}}^* \Gamma_{\rho^* {\mathcal F}} \cap \tilde{E}_{red}
= \pi_{\tilde{E}}^* \Gamma_{\rho^* {\mathcal L}} \cap
\tilde{E}_{red}
$$
for each irreducible component $\tilde{E}$ of
$\pi_{\tilde{C}}^{-1}(0)$, and from lemma~\ref{lem-polar-ram}, we
deduce that
$$
\pi_{\tilde{E}}^* \tilde{\Gamma}_{\mathcal F} \cap \tilde{E}_{red} =
\pi_{\tilde{E}}^* {\Gamma}_{\rho^* \mathcal F} \cap \tilde{E}_{red};
\ \ \pi_{\tilde{E}}^* \tilde{\Gamma}_{\mathcal L} \cap
\tilde{E}_{red} = \pi_{\tilde{E}}^* {\Gamma}_{\rho^* \mathcal L}
\cap \tilde{E}_{red}
$$
if $\tilde{E} \neq \tilde{E}_1$. Consequently,  $\pi_{\tilde{E}}^*
\tilde{\Gamma}_{\mathcal F} \cap \tilde{E}_{red} = \pi_{\tilde{E}}^*
\tilde{\Gamma}_{\mathcal L} \cap \tilde{E}_{red}$ provided that
$\tilde{E} \neq \tilde{E}_1$.

Moreover, the sets $\pi_{\tilde{E}_1}^* \tilde{\Gamma}_{\mathcal F}
\cap \tilde{E}_{1,red} \smallsetminus \pi_{\tilde{E}_1}^* \tilde{C}
\cap \tilde{E}_{1,red}$ and $\pi_{\tilde{E}_1}^*
\tilde{\Gamma}_{\mathcal F} \cap \tilde{E}_{1,red} \smallsetminus
\pi_{\tilde{E}_1}^* \tilde{C} \cap \tilde{E}_{1,red}$ have always
$b_{\tilde{E}_1}-1$ different points by lemma~\ref{lema-ptosE1}.
Then the result follows straightforward applying the criterion given
in proposition~\ref{criterio-adjunto-perfecto}.
\end{proof}
Now we are ready to prove proposition~\ref{prop:A-B}:
\begin{proof}[Proof of proposition~\ref{prop:A-B}]
Let ${\mathcal F}$ be a foliation in ${\mathbb G}_{C}^*$. By the
results of section~\ref{sec:non-sp}, it is clear that
$$
\epsilon(\Gamma_{\rho^{\ast}{\mathcal F}} \cap
\rho^{-1}C)=\chi_{\rho^{-1}C} \text{ if and only if, }
\lambda^*=\lambda(\rho^*{\mathcal F}) \in U_{\rho^{-1}C}=
\bigcap_{\tilde{E} \in B(\tilde{C})} U_{\tilde{C}}^{\tilde{E}},
$$
where $B(\tilde{C})$ is the set of bifurcation divisors of
$\pi_{\tilde{C}}^{-1}(0)$ and $U_{\tilde{C}}^{\tilde{E}} \subset
{\mathbb P}_{\mathbb C}^{m-1}$ are the Zariski-open sets defined in
section~\ref{sec:non-sp}. From lemmas~\ref{lema-ptosE1} and
\ref{lem-polar-ram} we deduce that
$$ \epsilon(\rho^{-1} \Gamma_{\mathcal F} \cap \rho^{-1}C)
= \chi_{\rho^{-1}C} \text{ if and only if } \lambda^* \in
\bigcap_{\tilde{E} \in B(\tilde{C}) \smallsetminus \{ \tilde{E}_1
\}} U_{\tilde{C}}^{\tilde{E}}.
$$
Consequently property $(A)$ implies $(B)$.

Assume now that $C$ has at most two different tangent lines, i.e.,
$b_{E_1}=b_{\tilde{E}_1} \leq 2$. If $b_{\tilde{E}_1}=1$, then
$\tilde{E}_1$ is not a bifurcation divisor. If $b_{\tilde{E}_1}=2$,
we can see that $U_{\tilde{C}}^{\tilde{E}_1}={\mathbb P}_{\mathbb
C}^{m_0(C)-1}$ (see its definition in section~\ref{sec:non-sp}). It
follows that $(A)$ and $(B)$ are equivalent when $C$ has at most two
different tangent lines.
\end{proof}
The set $U_C$ is equal to the set of  $\lambda$ such that each
${\mathcal F} \in {\mathbb G}_{C,\lambda}^{\ast}$ is a
Zariski-general foliation. It is an open subset of ${\mathbb
P}_{\mathbb C}^{r-1}$ but it could be empty. In fact, remark that
$\lambda=(\lambda_1,\ldots,\lambda_r) \in U_C$ if and only if,
$$
\lambda^* = (\overbrace{\lambda_{1}, \ldots,\lambda_{1}}^{n^1},
\ldots, \overbrace{\lambda_{r}, \ldots, \lambda_{r}}^{n^{r}}) \in
\bigcap_{\tilde{E} \in B(\tilde{C}) \atop  \tilde{E} \neq
\tilde{E}_1 } U_{\tilde{C}}^{\tilde{E}} \subset {\mathbb P}_{\mathbb
C}^{m_0(C)-1}
$$
where $n^i=m_0(C_i)$ for $i=1,\ldots,r$. The theorem
\ref{th:U-C-no-vacio} characterizes the equisingularity types
$\epsilon(C)$ such that $U_C \neq \emptyset$.

\begin{proof}[Proof of theorem~\ref{th:U-C-no-vacio}]
Let us see that,  for each bifurcation divisor $E$ of $G(C)$, we can
construct an open set $U_C^E \subset {\mathbb P}_{\mathbb C}^{r-1}$
such that
$$
\lambda \in U_C \text{ if and only if } \lambda \in \bigcap_{E \in
B(C)} U_{C}^E \text{ and } \sum_{i=1}^{r} k_i \lambda_i \neq 0
\text{ for } k \in R_{\epsilon(C)}.
$$
Moreover, we prove that a necessary and sufficient condition to
assure that each $U_C^E$ is non-empty is that $C$ has a kind
equisingularity type.

Consider a logarithmic foliation ${\mathcal L}_{\lambda} \in
{\mathbb G}_{C}^*$. Denote by $\Gamma_{\lambda}$ a generic polar
curve of ${\mathcal L}_{\lambda}$ and put $\tilde{\Gamma}_{\lambda}=
\rho^{-1} \Gamma_{\lambda}$. Take a bifurcation divisor $E$ of
$G(C)$ and let $\tilde{E}$ be any bifurcation divisor of
$G(\tilde{C})$ associated to $E$. Let us determine the conditions
over $\lambda$ which are equivalent to the fact that the set
$\pi_{\tilde{E}}^* \tilde{\Gamma}_{\lambda} \cap \tilde{E}_{red}
\smallsetminus \pi_{\tilde{E}}^* \tilde{C} \cap \tilde{E}_{red} $
has exactly $b_{\tilde{E}}-1$ different points. By
lemma~\ref{lema-ptosE1}, we only need to check this condition for
$\tilde{E} \neq \tilde{E}_1$ and hence, by
lemma~\ref{lem-polar-ram}, we have that
$$
\pi_{\tilde{E}}^* \tilde{\Gamma}_{\lambda} \cap \tilde{E}_{red}
\smallsetminus \pi_{\tilde{E}}^* \tilde{C} \cap \tilde{E}_{red} =
\pi_{\tilde{E}}^* {\Gamma}_{\lambda^*} \cap \tilde{E}_{red}
\smallsetminus \pi_{\tilde{E}}^* \tilde{C} \cap \tilde{E}_{red}
$$
where $\Gamma_{\lambda^*}$ is a generic polar curve of ${\mathcal
L}_{\lambda^*}= \rho^* {\mathcal L}_{\lambda}$.

Up to a coordinate change, we can assume that $(u,v)$ are
coordinates adapted to $\tilde{E}$. Let $\pi_{\tilde{E}} :
\tilde{M}_{\tilde{E}} \rightarrow ({\mathbb C}^2,0)$ be the morphism
reduction of $\pi_{\tilde{C}}$ to $\tilde{E}$ and take coordinates
$(u_p,v_p)$
 in the first chart of $\tilde{E}_{red}$ such that
$\tilde{E}_{red}=(u_p=0)$ and
$\pi_{\tilde{E}}(u_p,v_p)=(u_p,u_p^pv_p)$. Consider the 1-form
$$
\omega_{\lambda^*}^{\tilde{E}} = A_{\lambda^*}^{\tilde{E}} (u_p,v_p)
du_p + u_p B_{\lambda^*}^{\tilde{E}}(u_p,v_p) dv_p
$$
such that the strict transform $\pi_{\tilde{E}}^*{\mathcal
L}_{\lambda^*}$ is defined by $\omega_{\lambda^*}^{\tilde{E}}=0$. By
the results of section~\ref{sec:non-sp}, we know that the singular
points of $\pi_{\tilde{E}}^*{\mathcal L}_{\lambda^*}$ in the first
chart of $\tilde{E}_{red}$ are given by $u_p=0$ and
$A_{\lambda^*}^{\tilde{E}}(0,v_p)=0$ and the points of
$\pi_{\tilde{E}}^* \Gamma_{\lambda^*} \cap \tilde{E}_{red}$ are
given by $u_p=0$ and $B_{\lambda^*}^{\tilde{E}}(0,v_p)=0$. Denote by
$\{R_{1}^{\tilde{E}}, \ldots,R_{b_{\tilde{E}}}^{\tilde{E}}\}$  the
points of the set $\pi_{\tilde{E}}^*\tilde{C} \cap \tilde{E}_{red}$
with $R_{i}^{\tilde{E}}=(0,c_i^{\tilde{E}})$ in the coordinates
$(u_p,v_p)$. Note that these points are also the  singular points of
$\pi_{\tilde{E}}^* {\mathcal L}_{\lambda^*}$ in the first chart of
$\tilde{E}_{red}$. We deduce that, up to divide by a constant, we
have that
$$
A_{\lambda^*}^{\tilde{E}}(0,v)= \prod_{i=1}^{b_{\tilde{E}}}
(v-c_{i}^{\tilde{E}})^{r_i},
$$
where $r_i=m_{R_i^{\tilde{E}}}(\pi_{\tilde{E}}^* \tilde{C})$. We put
$A^{\tilde{E}}(v)=A_{\lambda^*}^{\tilde{E}}(0,v)$. Moreover, the
points of the set $\pi_{\tilde{E}}^* {\Gamma}_{\lambda^*} \cap
\tilde{E}_{red} \smallsetminus \pi_{\tilde{E}}^* \tilde{C} \cap
\tilde{E}_{red}$ are given by $u_p=0$ and
$H_{\lambda^*}^{\tilde{E}}(v_p)=0$ with
$$
H_{\lambda^*}^{\tilde{E}}(v)=
\frac{B_{\lambda^*}^{\tilde{E}}(0,v)}{\prod_{i=1}^{b_{\tilde{E}}}
(v-c_{i}^{\tilde{E}})^{r_i-1}}.
$$
The polynomial $H_{\lambda^*}^{\tilde{E}}(v)$ has degree
$b_{\tilde{E}}-1$ as a polynomial  in $v$ and its coefficients
depend linearly on $\lambda$;  we denote $H_{\lambda}^{\tilde{E}}(v)
= H_{\lambda^*}^{\tilde{E}}(v)$. Let $D^{\tilde{E}}(\lambda)$ be the
discriminant of $H_{\lambda}^{\tilde{E}}(v)$ as a polynomial in $v$
and we define  $U_{C}^E$ to be the set of $\lambda$ such that
$D^{\tilde{E}}(\lambda) \neq 0$ for all divisor $\tilde{E} \in
B(\tilde{C})$ associated to $E$. Let us show that each set $U_{C}^E$
is a non-empty Zariski open set if and only if $C$ has a kind
equisingularity type.

First we compute the polynomials above in terms of the Puiseux
series of the branches of $C$. The expression  of the polynomials
$A^{\tilde{E}}(v)$ and $B_{\lambda^*}^{\tilde{E}}(0,v)$ for a
logarithmic foliation with only non-singular separatrices in terms
of the parameterizations of its separatrices was described in the
proof of proposition~\ref{prop:U_C-log}. To compute these
polynomials in our situation we must take into account that the
curve $\tilde{C}$ is obtained by ramification from $C=\cup_{i=1}^r
C_i$. Consider a Puiseux series $y^i (x) = \sum_{s \geq n^i} a_{s}^i
x^{s/n^i}$ for each curve $C_i$ where $n^i=m_0(C_i)$. Thus all the
Puiseux series of $C_i$ are given by
$$
y_{j}^{i}(x) = \sum_{s \geq n^{i}} a_{s}^{i} (\varepsilon_i)^{sj}
x^{s/n^{i}}, \ \ \ \text{ for }j=1,2,\ldots, n^i,
$$
where $\varepsilon_i$ is a primitive $n^i$-root of the unity. Put
$v_{j}^i (u) = y_j^i(u^n)$. Then $\rho^{-1}C_i=
\{\sigma_{j}^{i}\}_{j=1}^{n^i}$ where
$\sigma_{j}^{i}=(v-v_{j}^{i}(u)=0)$.

Let $\{\tilde{E}^{l}\}_{i=1}^{\underline{n}_E}$ be the vertices of
$G(\tilde{C})$ associated to $E$ and assume that
$\tilde{E}=\tilde{E}^{l}$ for a certain $l \in \{1, \ldots,
\underline{n}_E\}$. By the results of
appendix~\ref{ap:ramificacion}, we know that the choice of a vertex
$\tilde{E}^l$ is equivalent to the choice of a $\underline{n}_E$-th
root $\xi_l$ of the unity. Given any $i \in I_E$, we denote
$e_{E}^{i}= n^{i}/\underline{n}_E$ and we consider
$\{\zeta_{ilt}\}_{t=1}^{e_{E}^{i}}$ the $e_{E}^{i}$-th roots of
$\xi_l$. Thus, if we denote by
$\{\sigma_{lt}^{i}\}_{t=1}^{e_{E}^{i}}$  the branches of
$\rho^{-1}C_i$ such that $\tilde{E}^l$ belongs to their geodesics,
then $\sigma_{lt}^{i}=(v-\eta_{lt}^{i}(u)=0)$ where
$$
\eta_{lt}^{i}(u)=\sum_{s \geq n^i} a_{s}^{i} (\zeta_{ilt})^s
u^{sn/n^i}, \text{ for } t=1, \ldots, e_{E}^{i}.
$$
The use of the expressions above to compute the polynomials
$A^{\tilde{E}^l}(v)$ and
$B_{\lambda}^{\tilde{E}^l}(v)=B_{\lambda^*}^{\tilde{E}^l}(0,v)$
gives that
\begin{align}
A^{\tilde{E}^l}(v)&= \prod_{i \in I_{E}} \prod_{t=1}^{e_{E}^i} (v -
a_{n^i v(E)}^i (\zeta_{ilt})^{n^i v(E)}) \tag{$*_1$}
\label{eq:A-ram}
\\
B_{\lambda}^{\tilde{E}^l}(v)&= \sum_{i \in I_E} \lambda_i
\prod_{j\in I_E \atop j \neq i} \prod_{t=1}^{e_E^j} (v-a_{n^j
v(E)}^j \zeta_{jlt}^{n^j v(E)})\sum_{t=1}^{e_{E}^i} \prod_{k=1 \atop
k \neq t}^{e_E^i} (v-a_{n^i v(E)}^{i} \zeta_{ilk}^{n^i v(E)})
\tag{$*_2$} \label{eq:B-ram}
\end{align}
Since both polynomials only depend on the invariants associated to
$E$, we consider the three possibilities for a divisor $E$ of $G(C)$
in order to obtain a more precisely expression of them:

\hspace{.5cm}  (i) {\em  $E$ is a contact divisor:\/} we have that
$v(E) = m_E /\underline{n}_E$ and $n_E=1$. Then $n^i v(E)= e_{E}^i
m_E$ for each $i \in I_E$ and consequently
$(\zeta_{ilt})^{n^{i}v(E)} = \xi_{l}^{m_{E}}$ for each $t \in  \{1,
\ldots, e_{E}^i \}$. Thus we have that
\begin{align*}
A^{\tilde{E}^l}(v) &= \prod_{i \in I_E} (v-a_{n^i v(E)}^{i}
\xi_{l}^{m_E})^{e_{E}^i} \\
B_{\lambda}^{\tilde{E}^l}(v) &=  \prod_{j \in I_{E}} (v-
a_{n^jv(E)}^{j} \xi_{l}^{m_E})^{e_{E}^{j}-1} \sum_{i \in I_E}
\lambda_i e_{E}^i \prod_{j \in I_{E} \atop j \neq i} (v-
a_{n^jv(E)}^{j} \xi_{l}^{m_E})
\end{align*}
Denote by $I_{\tilde{E}^l}^s= \{ i \in I_{E} \ : \ a_{n^i v(E)}^{i}
\xi_{l}^{m_E} = c_s^{\tilde{E}^l} \}$ for $s=1, \ldots,
b_{\tilde{E}^l}$. Thus $r_s= \sum_{i \in I_{\tilde{E}^l}^{s}}
e_{E}^i$ and we have that
$$
H_{\lambda}^{\tilde{E}^l}(v) = \sum_{i=1}^{b_{\tilde{E}^l}} (
\sum_{s \in I_{\tilde{E}^l}^{i}} \lambda_s e_{E}^{s}) \prod_{j=1
\atop j\neq i}^{b_{\tilde{E}^l}} (v-c_{j}^{\tilde{E}^l})
$$
which is a polynomial of degree $b_{\tilde{E}^l}-1$ in $v$. Observe
that $b_{\tilde{E}^l} = b_{E}$. The discriminant
$D^{\tilde{E}^{l}}(\lambda)$  of $H_{\lambda}^{\tilde{E}^l}(v)$ as a
polynomial in $v$ is a non-zero polynomial. Hence, the set
$U_{C}^{E}= \{\lambda \ : \  D^{\tilde{E}^l}(\lambda) \neq 0 \text{
for } l=1, \ldots, \underline{n}_E\}$ is a non-empty Zariski open
set.

\vspace{\baselineskip} \hspace{.5cm}(ii) {\em $E$ is a Puiseux
divisor with a dead arc:\/} we have that $v(E)=m_E /\underline{n}_E
n_E$ with $n_E>1$ and $(m_{k_E+1}^{i},n_{k_{E}+1}^{i})=(m_E,n_E)$
for each $i \in I_E$. It follows that $n^i v(E) = e_E^i m_E /n_E$
and the set $\{\zeta_{ilt}^{n^iv(E)}\}_{t=1}^{e_{E}^i}$ has $n_E$
different values which coincide with the $n_E$-th roots
$\{\theta_{lt}\}_{t=1}^{n_E}$ of $\xi_{l}^{m_E}$. Moreover, we have
that
$$\prod_{s=1}^{n_E} (v-a_{n^iv(E)}^i \theta_{ls}) =
v^{n_E}-\alpha_{\tilde{E}^l}^{i} \text{ with }
\alpha_{\tilde{E}^l}^{i} = (a_{n^iv(E)}^{i})^{n_E} \xi_{l}^{m_E}
$$
and $\sum_{t=1}^{n_E} \prod_{p=1 \atop p \neq
t}^{n_E}(v-a_{n^{i}v(E)}^{i} \theta_{lp}) = n_E v^{n_E-1}$.  Thus
the expressions \eqref{eq:A-ram} and \eqref{eq:B-ram} become
\begin{align*}
A^{\tilde{E}^l}(v) & =
\prod_{i \in I_E} (v^{n_E}-\alpha_{\tilde{E}^l}^{i})^{e_{E}^{i}/n_E} \\
B_{\lambda}^{\tilde{E}^l}(v) &= n_E v^{n_E-1} \prod_{i \in I_E}
(v^{n_E}-\alpha_{\tilde{E}^l}^{i})^{\frac{e_{E}^{i}}{n_E}-1} \sum_{i
\in I_E} \lambda_i \frac{e_{E}^{i}}{n_E} \prod_{j \in I_E \atop j
\neq i} (v^{n_E}-\alpha_{\tilde{E}^l}^{j})
\end{align*}
In this case we have that $b_{\tilde{E}^l}= n_E(b_E-1)$ and hence
there are exactly $b_E-1$ different values
$\{\phi^{\tilde{E}_l}_s\}_{s=1}^{b_E-1}$ in the set
$\{\alpha_{\tilde{E}_l}^{i} \}_{i \in I_E}$. Denote
$I_{\tilde{E}^l}^{s} = \{ i \in I_E :  \ \alpha_{\tilde{E}^l}^{i} =
\phi_s^{\tilde{E}^l}\}$ and $r_s= \sum_{i \in I_{\tilde{E}^l}^{s}}
e_{E}^{i}/n_E$. Then we have that $A^{\tilde{E}^l}(v)=
\prod_{s=1}^{b_E-1} (v^{n_E}-\phi_{s}^{\tilde{E}^l})^{r_s}$ and
$$
H_{\lambda}^{\tilde{E}^l}(v) = v^{n_E-1} \sum_{s=1}^{b_E-1} (
\sum_{i \in I_{\tilde{E}^l}^{s}} \lambda_i e_{E}^{i} ) \prod_{j=1
\atop j\neq s}^{b_E-1} (v^{n_E}-\phi_{j}^{\tilde{E}^l}).
$$
In this situation $D^{\tilde{E}^l}(\lambda) \not \equiv 0$ if and
only if $n_E =2$. Hence, we conclude that $U_C^E$ is a non-empty
Zariski open set if and only if  $C$ has a kind  equisingularity
type.

 \vspace{\baselineskip} \hspace{.5cm}(iii) {\em $E$ is a
Puiseux divisor without a dead arc:\/} we have that $v(E)=m_E
/\underline{n}_E n_E$ with $n_E >1$ and $b_{\tilde{E}^l}= 1+n_E
(b_{E}-1)$. We know that $(m_E,n_E)=(m_{k_E+1}^{i},n_{k_E+1}^{i})$
for each $i \in I_{E}^{\ast}$ and $a_{n^{i}v(E)}^j=0$ for $i \in I_E
\smallsetminus I_{E}^*$ (see appendix~\ref{ap:grafodual}). Denote by
$r_0= \sharp(I_E \smallsetminus I_E^*)$. With similar arguments and
notations as in case (ii), we get that
\begin{align*}
A^{\tilde{E}^l}(v) & = v^{r_0} \prod_{i \in I_{E}^*} (v^{n_E} -
\alpha^{i}_{\tilde{E}^l})^{e_{E}^{i}/n_E} \\
B_{\lambda}^{\tilde{E}^l}(v) &= v^{r_0-1} \prod_{i \in I_{E}^*}
(v^{n_E}- \alpha^{i}_{\tilde{E}^l})^{\frac{e_{E}^{i}}{n_E}-1} &
\left\{ v^{n_E} \sum_{i \in I_{E}^*}  \lambda_i e_{E}^i \prod_{j
\in I_{E}^* \atop j \neq i} (v^{n_E}-\alpha_{\tilde{E}^l}^{j}) + \right.  \\
& & + \left. \prod_{j \in I_{E}^*}
(v^{n_E}-\alpha_{\tilde{E}^l}^{j}) (\sum_{i \in I_E \smallsetminus
I_{E}^*} \lambda_i) \right\}
\end{align*}
Let $\{\phi^{\tilde{E}_l}_{s}\}_{s=1}^{b_E-1}$ be the $b_{E}-1$
different values in the set $\{ \alpha_{i}^{\tilde{E}^l} \}_{i \in
I_{E}^*}$. Denote $I_{\tilde{E}^l}^{s} = \{ i \in I_E^* :  \
\alpha_{\tilde{E}^l}^{i} = \phi_s^{\tilde{E}^l}\}$  and $r_s=
\sum_{i \in I_{\tilde{E}^l}^{s}} e_{E}^{i}/n_E$. Thus we have that
\begin{align*}
A^{\tilde{E}^l}(v) &= v^{r_0} \prod_{i=1}^{b_E-1}
(v^{n_E}-\phi_s^{\tilde{E}^l})^{r_i} \\
H_{\lambda}^{\tilde{E}^l}(v) &= v^{n_E} \sum_{s=1}^{b_E-1} (\sum_{i
\in I_{\tilde{E}^l}^{s}} \lambda_i e_{E}^i) \prod_{j=1 \atop j \neq
s}^{b_E-1} (v^{n_E}-\phi^{\tilde{E}^l}_j) + (\sum_{j \in I_E
\smallsetminus I_{E}^*} \lambda_j) \prod_{s=1}^{b_E-1}
(v^{n_E}-\phi_s^{\tilde{E}^l}).
\end{align*}
It is clear that in this case $D^{\tilde{E}^l}(\lambda) \not \equiv
0$ for each $l=1,\ldots,\underline{n}_E$. Consequently, $U_C^E$ is a
non-empty Zariski open set.

We conclude that a necessary and sufficient condition to assure that
all the sets $U_C^E$ are non-empty Zariski open sets is that  $C$
has a kind equisingularity type and  the result follows
straightforward.
\end{proof}
With similar arguments to the ones in the proof above  we can show
that:
\begin{Corollary}
The following statements are equivalent:
\begin{itemize}
    \item The curve $C$ has a kind equisingularity type;
    \item There exists a foliation ${\mathcal F} \in {\mathbb
G}_{C}^*$ such that $\rho^\ast {\mathcal F}$ is Zariski-general.
\end{itemize}
\end{Corollary}
In particular, if ${\mathcal F} \in {\mathbb G}_{C,\lambda}^*$ with
$\lambda \in U_C$, the equisingularity type of a generic polar curve
$\Gamma_{\mathcal F}$ is completely determined in terms of $C$ and
$\pi_C$ gives a reduction of singularities of $C \cup
\Gamma_{\mathcal F}$. Moreover, we get that the irreducible
components of $\Gamma_{\mathcal F}$ cut the exceptional divisor
$\pi_{C}^{-1}(0)$ as described in
proposition~\ref{prop:equis-adjperf}; we get a more specific
description than the one of  L\^e-Michel-Weber in \cite{Le-M-W}.

Observe that the property ``$\pi_C$ gives a reduction of
singularities of $\Gamma_{\mathcal F} \cup C$" does not imply that
$\mathcal F$ is a Zariski-general foliation. Moreover, this property
does not determine the equisingularity type of $\Gamma_{\mathcal F}
\cup C$ even if we fix $\lambda$.
\begin{Example}
Consider the foliations ${\mathcal F}_1$, ${\mathcal F}_2$ and
${\mathcal F}_3$ given by $\omega_i=0$ with
\begin{align*}
\omega_1& = -11 x^{10}dx+5y^4dy; \\
\omega_2&=11(-x^{10}+y^2x^6)dx+5(y^4-x^7y)dy;\\
\omega_3&=11(-x^{10}+yx^8)dx+5(y^4-x^9)dy
\end{align*}
respectively. All the foliations   have the same separatrix
$C=(y^5-x^{11}=0)$ which does not have a kind type of
equisingularity, therefore  ${\mathcal F}_1$, ${\mathcal F}_2$ and
${\mathcal F}_3$ cannot be Zariski-general foliations. The generic
polar curves $\Gamma_{{\mathcal F}_1}$, $\Gamma_{{\mathcal F}_2}$
and $\Gamma_{{\mathcal F}_3}$ are not equisingular but the minimal
reduction of singularities of $C$ is also a reduction of
singularities of the curves $\Gamma_{{\mathcal F}_1}$,
$\Gamma_{{\mathcal F}_2}$ and $\Gamma_{{\mathcal F}_3}$.

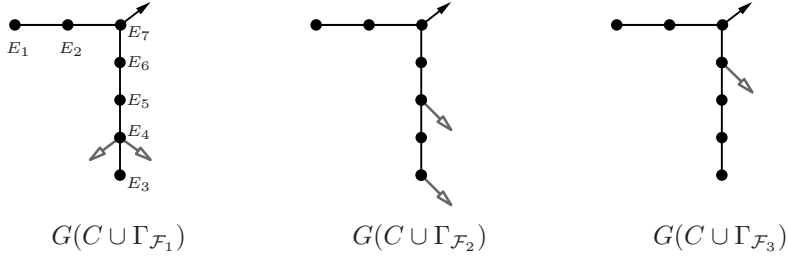
\vspace{\baselineskip}

\begin{center}
\begin{texdraw}
\arrowheadtype t:F \arrowheadsize l:0.08 w:0.04 \drawdim mm \setgray
0

\move (10 40) \fcir f:0 r:.8 \rlvec (7 0) \fcir f:0 r:.8 \rlvec (7
0) \fcir f:0 r:.8 \avec (28 43) \move(24 40) \rlvec (0 -5) \fcir f:0
r:.8 \rlvec (0 -5) \fcir f:0 r:.8 \rlvec (0 -5) \fcir f:0 r:.8
\rlvec (0 -5) \fcir f:0 r:.8

\move(9 36) \htext{\tiny{$E_1$}} \move(16 36) \htext{\tiny{$E_2$}}
\move(25 18) \htext{\tiny{$E_3$}} \move(25 25) \htext{\tiny{$E_4$}}
\move(25 29) \htext{\tiny{$E_5$}} \move(25 34) \htext{\tiny{$E_6$}}
\move(25 38) \htext{\tiny{$E_7$}}

\move(50 40)
\fcir f:0 r:.8 \rlvec (7 0) \fcir f:0 r:.8 \rlvec (7 0) \fcir f:0
r:.8 \avec (68 43) \move(64 40) \rlvec (0 -5) \fcir f:0 r:.8 \rlvec
(0 -5) \fcir f:0 r:.8 \rlvec (0 -5) \fcir f:0 r:.8 \rlvec (0 -5)
\fcir f:0 r:.8

\move(90 40)
\fcir f:0 r:.8 \rlvec (7 0) \fcir f:0 r:.8 \rlvec (7 0) \fcir f:0
r:.8 \avec (108 43) \move(104 40) \rlvec (0 -5) \fcir f:0 r:.8
\rlvec (0 -5) \fcir f:0 r:.8 \rlvec (0 -5) \fcir f:0 r:.8 \rlvec (0
-5) \fcir f:0 r:.8

\linewd 0.3 \setgray 0.4 \arrowheadtype t:W \move(24 25) \avec(28
22) \move(24 25) \avec(20 22)

\move(64 30) \avec(68 26) \move(64 20) \avec(68 16)

\move(104 35) \avec(108 31)

\move(24 25) \fcir f:0 r:.8 \move(64 30) \fcir f:0 r:.8 \move(64 20)
\fcir f:0 r:.8 \move(104 35) \fcir f:0 r:.8

\move(15 10) \htext{$G(C \cup \Gamma_{{\mathcal F}_1})$} \move(55
10) \htext{$G(C \cup \Gamma_{{\mathcal F}_2})$} \move(95 10)
\htext{$G(C \cup \Gamma_{{\mathcal F}_3})$}

\end{texdraw}
\end{center}
\end{Example}

\appendix
\section{Equisingularity data: the dual graph}\label{ap:grafodual}
Let us recall the construction of the dual graph which is one of the
different ways to represent the equisingularity data of a plane
curve (see \cite{Bri-K} for more details). Let $C \subset ({\mathbb
C}^2,0)$ be a plane curve and $\pi_{C} : M \rightarrow ({\mathbb
C}^2,0)$ be its minimal reduction of singularities. The {\em dual
graph\/} $G(C)$ is constructed as follows:  each irreducible
component $E$ of $\pi_{C}^{-1}(0)$ is represented by a vertex which
we also call $E$ (we identify a divisor and its  associated vertex
in the dual graph). Two vertices are joined by an edge if and only
if the associated divisors intersect. Each irreducible component of
$C$ is represented by an arrow joined to the only divisor which
meets the strict transform of $C$ by $\pi_C$. If we give a weight to
each vertex $E$ of $G(C)$ equal to the  self-intersection of the
divisor $E \subset M$, this weighted dual graph is equivalent to the
equisingularity data of $C$.

We  denote by $E_1$ the irreducible component of $\pi_{C}^{-1}(0)$
obtained by the blowing-up of the origin. Thus the first divisor
$E_1$ gives an orientation to the graph $G(C)$. The {\em geodesic\/}
of a divisor $E$ is the path which joins the first divisor $E_1$
with the divisor $E$. The geodesic of a curve is the geodesic of the
divisor that meets the transform strict of the curve. In this way,
there is a partial order in the set of vertices of $G(C)$ given by
$E < E'$ if and only if the geodesic of $E'$ goes through $E$.

Let us introduce some notations concerning the dual graph of a
curve. Given a vertex $E$ of $G(C)$ we define the number $b_E$ as
follows: $b_E+1$ is the valence of $E$ if $E \neq E_1$ and $b_{E_1}$
is the valence of $E_1$. Observe that $b_{E_1}$  is the number of
different lines in the tangent cone of $C$. We say that $E$ is a
{\em bifurcation divisor\/} if $b_E \geq 2$ and $E$ is a {\em
terminal divisor\/} if $b_E=0$. A {\em dead arc\/} in $G(C)$ is an
arc which joins a bifurcation divisor with a terminal one without
passing through other bifurcation divisors. Observe that a
bifurcation divisor can belong only to one dead arc.

A {\em curvette\/} $\tilde{\gamma}$ of a divisor $E$ is a
non-singular curve transversal to $E$ at a non-singular point of
$\pi_{C}^{-1}(0)$. The projection $\gamma=\pi_{C}(\tilde{\gamma})$
is a germ of plane curve in $({\mathbb C}^2,0)$ and we say that
$\gamma$ is an $E$-curvette. We denote by $m(E)$ the multiplicity at
the origin of any $E$-curvette. Take
$\tilde{\gamma},\tilde{\gamma}'$ two curvettes of $E$ which
intersect $E$ in two different points, we denote by $v(E)$ the
coincidence ${\mathcal
C}(\pi_C(\tilde{\gamma}),\pi_C(\tilde{\gamma}'))$; then $v(E) <
v(E')$ if $E < E'$. Recall  that the {\em coin\-cidence\/}
${\mathcal C}(\gamma,\delta)$ between two irreducible curves
$\gamma$ and $\delta$ is defined as
$$
{\mathcal C}(\gamma,\delta) = \sup_{\ 1 \leq i \leq m_0(\gamma)
\atop 1 \leq j \leq m_0(\delta)} \{
ord_{x}(y^{\gamma}_i(x)-y^{\delta}_j(x))\ \}
$$
where $\{y_{i}^{\gamma}(x)\}_{i=1}^{m_0(\gamma)}$,
$\{y_{j}^{\delta}(x)\}_{j=1}^{m_0(\delta)}$ are the Puiseux series
of $\gamma$ and $\delta$ respectively.

Given any irreducible component $E$ of the exceptional divisor
$\pi^{-1}_C(0)$, we denote by $\pi_E:M_E \rightarrow ({\mathbb
C}^2,0)$ the {\em reduction of $\pi_C$ to $E$\/}, that is, the
morphism which satisfies that
\begin{itemize}
    \item there is a factorization $\pi_C= \pi_{E}' \circ \pi_E$
    where $\pi_{E}'$ and $\pi_E$ are composition of punctual
    blow-ups;
    \item the divisor $E$ is the strict transform by $\pi_{E}'$ of
    an irreducible component $E_{red}$ of $\pi_{E}^{-1}(0)$ and
    $E_{red} \subset M_E$ is the only component of $\pi_{E}^{-1}(0)$
    with self-intersection equal to $-1$.
\end{itemize}
It is clear that $\pi_E$ is obtained from $\pi_C$ by blowing-down
successively the divisors which are different from $E$ and whose
self-intersection is equal to $-1$. Take any curvette
$\tilde{\gamma}_E$ of $E$, then $\pi_{E}'(\tilde{\gamma}_E)$ is also
a curvette of $E_{red} \subset M_E$. Let $\{ \beta_0^E,\beta_1^E,
\ldots, \beta_{g(E)}^{E}\}$ be the characteristic exponents of
$\gamma_E = \pi_{C}(\tilde{\gamma}_E)$. It is clear that
$m(E)=\beta_0^E=m_0(\gamma_{E})$ and there are two possibilities for
the value $v(E)$:
\begin{itemize}
    \item[1.] either $\pi_E$ is the minimal reduction of
    singularities of $\gamma_E$ and then
    $v(E)=\beta_{g(E)}^E/\beta_0^E$. We say that $E$ is a {\em
    Puiseux divisor\/} for $\pi_C$.
    \item[2.] or $\pi_E$ is obtained by blowing-up $q \geq 1$ times
    after the minimal reduction of singularities of $\gamma_E$ and
    in this situation $v(E) = (\beta_{g(E)}^{E}+q \beta_0^E)
    /\beta_0^E$. We say that $E$ is a {\em contact divisor\/} for
    $\pi_C$.
\end{itemize}
Observe that $m(E)=m(E_{red})$ and $v(E)=v(E_{red})$. Moreover, $E$
can belong to a dead arc only if it is a Puiseux divisor.

Consider a bifurcation divisor $E$ of $G(C)$ and let
$\{(m_1^E,n_1^E), (m_2^E,n_2^E), \ldots,
\newline (m_{g(E)}^E,n_{g(E)}^E)\}$ be the Puiseux pairs of an $E$-curvette
$\gamma_E$, we denote
$$n_E =
\left\{
  \begin{array}{ll}
    n_{g(E)}, & \hbox{ if $E$ is a Puiseux divisor;} \\
    1, & \hbox{otherwise,}
  \end{array}
\right.
$$
and $\underline{n}_E=m(E)/n_E$. Observe that, if $E$ belongs to a
dead arc with terminal divisor $F$, then $m(F)=\underline{n}_E$. We
define $k_E$ to be
$$
k_E =\left\{
        \begin{array}{ll}
          g(E)-1, & \hbox{if $E$ is a Puiseux divisor;} \\
          g(E), & \hbox{if $E$ is a contact divisor.}
        \end{array}
      \right.
$$

Let us explain these notations in terms of the equisingularity data
of the curve $C= \cup_{i=1}^{r}C_i$. Denote by
$\{(m_l^i,n_l^i)\}_{l=1}^{g_i}$ the Puiseux pairs of $C_i$ and by
$\{\beta_0^i,\beta_1^i,\ldots,\beta_{g_i}^{i}\}$ its characteristic
exponents. Denote $I=\{1,2,\ldots,r\}$ and let $I_E$ be the set of
indices $i \in I$ such that $E$ belong to the geodesic of $C_i$.
Take $i \in I_E$.  There are several possibilities for the value of
$v(E)$ depending on $E$:
\begin{itemize}
    \item[(i)]  If $E$ is a contact divisor, then there exists $j \in I_E$
    such that $v(E)={\mathcal C}(C_i,C_j)$.
    \item[(ii)] If $E$ is a Puiseux divisor which  belongs to a dead arc, then
    $v(E)=\beta_{k_E+1}^i/\beta_0^i$.
    \item[(iii)] If $E$ is a Puiseux divisor which does not belong to a dead arc,
    we denote by $I_E^*$ the set of indices $i \in I_E$ such that
    $v(E)=\beta_{k_E+1}^i/\beta_0^i$. Then ${\mathcal C}(C_i,C_j) =
    v(E)$ for $i \in I_E^*$ and $j \in I_E \smallsetminus I_E^*$.
    Moreover, ${\mathcal C}(C_j,C_l) > v(E)$ if $j,l \in I_E \smallsetminus
    I_{E}^*$.
\end{itemize}
Consequently,  we have that $(m_l^{i},n_l^{i})=(m_{l}^E,n_{l}^E)$,
for $l=1, \ldots, k_E$, and $\underline{n}_{E}= n_1^i \cdots
n_{k_E}^i$ for any $i \in I_E$.

\section{Ramification}\label{ap:ramificacion}
Consider a plane curve $C=\cup_{i=1}^r C_i \subset ({\mathbb
C}^2,0)$. Let $\rho: ({\mathbb C}^2,0) \rightarrow ({\mathbb
C}^2,0)$ be any $C$-ramification, that is, $\rho$ is transversal to
C and  $\tilde{C} = \rho^{-1}C$ has only non-singular irreducible
components. Assume that the ramification is given by $x=u^n, y=v$.

Denote by $\{(m_l^i,n_l^i)\}_{l=1}^{g_i}$ the Puiseux pairs of $C_i$
and by $\{\beta_0^i,\beta_1^i,\ldots,\beta_{g_i}^{i}\}$ the
characteristic exponents of $C_i$. If $n^{i}=m_0(C_i)$, then it is
necessary that $n \equiv 0 \mod(n^1,n^2,\ldots,n^r)$ in order to
have that $\tilde{C}$ has only non-singular irreducible components.
Moreover, the number of irreducible components of $\tilde{C}$ is
equal to $m_0(C)= n^1+\cdots+n^r$. More precisely, each curve
$\rho^{-1}C_i$ has exactly $n^i$ irreducible components. In fact,
let $y^i(x) = \sum_{l \geq n^i} a_{l}^{i} x^{l/n^{i}}$ be a Puiseux
series of $C_i$, thus all its Puiseux series are given by
$$
y_{j}^{i}(x) = \sum_{l \geq n^{i}} a_{l}^{i} \varepsilon_i^{lj}
x^{l/n^{i}} \ \ \ \text{ for }j=1,2,\ldots, n^i,
$$
where $\varepsilon_i$ is a primitive $n^i$-root of the unity. Then
$f_i(x,y)=\prod_{l=1}^{n^{i}} (y-y_{l}^{i}(x))$ is a reduced
equation of $C_i$. If we put $v_j^{i}(u)=y_{j}^{i}(u^n)$, then
$v_{j}^{i}(u) \in {\mathbb C}\{u\}$ since $n/n^{i} \in {\mathbb N}$.
It is clear that the curve $\sigma_{j}^{i}=(v-v_j^{i}(u)=0)$ is
non-singular and it is one of the irreducible components of
$\rho^{-1}C_i$. Then
$$
g_i(u,v)= f_i(u^n,v)=\prod_{l=1}^{n^{i}} (v-v_{l}^{i}(u))
$$
is an equation of $\rho^{-1}C_i$. We conclude that the irreducible
components $\{\sigma_{j}^{i}\}_{j=1}^{n^{i}}$ of $\rho^{-1}C_i$ are
in bijection with the Puiseux series of $C_i$.

It is well-known that the equisingularity type of a curve $C$ is
determined by the characteristic exponents
$\{\beta_0^{i},\beta_1^{i},\ldots,\beta_{g_i}^{i}\}_{i=1}^{r}$ of
its irreducible components and the intersection multiplicities
$\{(C_i,C_j)_0\}_{i\neq j}$. Let us show that we can obtain all this
information from $\rho^{-1}C$.  The next lemma states the
relationship between the intersection multiplicity
$(\gamma,\delta)_0$ and the coincidence ${\mathcal
C}(\gamma,\delta)$ (see Zariski \cite{Zar}, prop. 6.1 or Merle
\cite{Mer}, prop. 2.4):
\begin{Lemma}\label{lema:coinc-multinterseccion}
Let $\gamma$ and $\delta$ be two germs of irreducible plane curves
of $({\mathbb C}^2,0)$. If $\{\beta_0,\beta_1,\ldots,\beta_g\}$ are
the characteristic exponents of $\gamma$ and $\alpha$ is a rational
number such that $\beta_q \leq \alpha < \beta_{q+1}$
($\beta_{g+1}=\infty$), then the following statements are
equivalent:
\begin{align*}
    1.& \ \ {\mathcal C}(\gamma,\delta) =
    \frac{\alpha}{m_0(\gamma)} \hspace{8cm}\\
    2. & \ \ \frac{(\gamma,\delta)_0}{m_0(\delta)}
    =\frac{\bar{\beta}_q}{n_1 \cdots n_{q-1}} +
    \frac{\alpha-\beta_q}{n_1 \cdots n_q}
\end{align*}
where $\{(m_i,n_i)\}_{i=1}^{g}$ are the Puiseux pairs of $\gamma$
($n_0=1$) and $\{\bar{\beta}_0,\bar{\beta}_1,\ldots,\bar{\beta}_q\}$
is a minimal system of generators of the semigroup $S(\gamma)$ of
$\gamma$.
\end{Lemma}
In particular, the equisingularity type of $C$ is also determined by
the characteristic exponents of each $C_i$ and the coincidences
$\{{\mathcal C}(C_i,C_j)\}_{i \neq j}$. Let us show that these data
could be obtained from $\rho^{-1}C$. Given an irreducible component
$\sigma$ of $\rho^{-1}C$, we take an equation $(v-v^{\sigma}(u)=0)$
of $\sigma$ with $v^{\sigma}(u) = \sum_{l \geq 1} a_l^{\sigma} u^{l}
\in {\mathbb C} \{u\}$. Given two irreducible components
$\sigma,\sigma'$ of $\rho^{-1}C$, we say that they are equivalent
$\sigma \sim \sigma'$ if and only if $(a_j^{\sigma})^n =
(a_{j}^{\sigma'})^n$ for all $j \in {\mathbb N}$. Denote by
$[\sigma]$ the equivalence classes of a curve $\sigma$. Thus the
number of irreducible components $r$ of $C$ is equal to the number
of equivalence classes for the irreducible components of
$\rho^{-1}C$. Let  $[\sigma^1],\ldots, [\sigma^r]$ be these
equivalence classes. Up to reorder, we can assume that
$[\sigma^{i}]$ corresponds to $\rho^{-1}C_i$, for $i=1,\ldots, r$.
Thus the multiplicity $n^i$ of $\rho^{-1}C_i$ is equal to the number
of elements in the equivalence class $[\sigma^{i}]$. We put
$\rho^{-1}C_i= \{\sigma_{l}^{i}\}_{l=1}^{n^{i}}$. Hence $\beta_0^i=
n^i$ and the other characteristic exponents of $C_i$ are obtained
from the computation of the coincidences among the curves in the
equivalence class $[\sigma^{i}]$ since
$$
\{ {\mathcal C}(\sigma_j^{i},\sigma_{l}^{i}) \ : \ j \neq l \}=
\{\beta_1^{i},\ldots,\beta_{g_i}^{i}\}.
$$
Thus we only need to compute the coincidences between any two
branches $C_i$ and $C_j$. But they are obtained from the following
equality
\begin{equation}\label{eq:coincidencias-ram}
{\mathcal C}(C_i,C_j) = \frac{1}{n} \sup_{1 \leq l \leq n^{i} \atop
1 \leq s \leq n^{j} } \{ {\mathcal
C}(\sigma_{l}^{i},\sigma_{s}^{j})\},
\end{equation}
which is true for any two irreducible curves. Hence we conclude that
the equisingularity data of $C$ can be recovered from $\rho^{-1}C$.

\vspace{\baselineskip}
\paragraph{\sc Ramification of the dual graph}
Let $\pi_{C}:M \rightarrow ({\mathbb C}^2,0)$ be the minimal
reduction of singularities of $C$ and denote by  $\pi_{\tilde{C}}:
\tilde{M} \rightarrow ({\mathbb C}^2,0)$ the minimal reduction of
singularities of $\tilde{C}=\rho^{-1}C$. Let us explain the
relationship between $G(C)$ and $G(\tilde{C})$.

Let $K_i$ be the geodesic in $G(C)$ of a branch $C_i$ of $C$ and let
$\tilde{K}_i$ be the sub-graph of $G(\tilde{C})$ corresponding to
the geodesics of the irreducible components
$\{\sigma_{l}^{i}\}_{l=1}^{n^i}$ of $\rho^{-1}C_i$. Let us see how
to construct $\tilde{K}_i$ from $K_i$. Observe first that, if
$\tilde{E}$ and $\tilde{E}'$ are two consecutive vertices of
$G(\tilde{C})$ with $\tilde{E} < \tilde{E}'$, then
$v(\tilde{E}')=v(\tilde{E})+1$. Thus, $G(\tilde{C})$ is completely
determined once we know  the bifurcation divisors, the order
relations  among them and the number of edges   which leave from
each bifurcation divisor. Denote by $B(\tilde{K}_i)$ and $B(K_i)$
the bifurcation vertices of $\tilde{K}_i$ and $K_i$ respectively. We
say that a vertex $\tilde{E}$ of $B(\tilde{K}_i)$ is associated to a
vertex $E$ of $B(K_i)$ if $v(\tilde{E})=n v(E)$.

Let $E$ be a vertex of $B(K_i)$. Assume first that $E$ is the first
bifurcation divisor of $B(K_i)$ and take $E'$ its consecutive vertex
in $B(K_i)$. Then $E$ has only one associated vertex $\tilde{E}$ in
$B(\tilde{K}_i)$ and there are two possibilities for the number of
edges which leave from it:
\begin{itemize}
    \item If $E$ is a Puiseux divisor,
    then there are $n_1^{i}$ edges which leave from $\tilde{E}$ in
    $\tilde{K}_i$; then  $E'$ has $n_1^i$ associated
    vertices in $B(\tilde{K}_i)$.
    \item If $E$ is a contact divisor, then
    there is only one edge which leave from $\tilde{E}$ in $\tilde{K}_i$ and thus
    $E'$ has only one vertex associated in $B(\tilde{K}_i)$.
\end{itemize}
Take now any vertex $E$ of $B(K_i)$ and assume that we know the part
of $\tilde{K}_i$ corresponding to the vertices of $K_i$ with
valuation $\leq v(E)$. Then there are $\underline{n}_E= n_1^i \cdots
n_{k_E}^i$ vertices $\{ \tilde{E}^{l} \}_{l=1}^{\underline{n}_E}$
associated to $E$ and
\begin{itemize}
    \item If $E$ is a Puiseux divisor,
    then there are $n_{k_E+1}$ edges which leave from each
    vertex $\tilde{E}_l$ in $\tilde{K}_i$.
    \item If $E$ is a contact divisor,
the there is only one edge which leaves from each vertex $\bar{E}_l$
in $\tilde{K}_i$.
\end{itemize}
The dual graph $G(\tilde{C})$ is constructed in the natural way by
gluing the graphs $\tilde{K}_i$. From the construction described
above, we deduced that
\begin{equation*}
b_{\tilde{E}} = \left\{
                  \begin{array}{ll}
                    b_E, & \hbox{if $E$ is a contact divisor;} \\
                    (b_E-1)n_E, & \hbox{if $E$ is a bifurcation divisor which
                                belong} \\
                                & \hbox{to a dead arc;} \\
                    (b_E-1)n_E +1 , & \hbox{if $E$ is a bifurcarion divisor which does
                                        not} \\
                                 & \hbox{belong to a dead arc.}
                  \end{array}
                \right.
\end{equation*}

Observe that, in general, non-bifurcation divisors  of $G(C)$ have
no associated divisors in $G(\tilde{C})$. Let us illustrate with
some examples the relationship between $G(C)$ and $G(\tilde{C})$:
\begin{Example}
Consider the curve $C=(y^2-x^3=0)$ and the rami\-fication
$\rho(u,v)=(u^2,v)$. Then $\tilde{C}$ has two irreducible components
given by $v-u^3=0$ and $v+u^3=0$. The next figure represents the
dual graphs of $C$ and $\rho^{-1}C$:
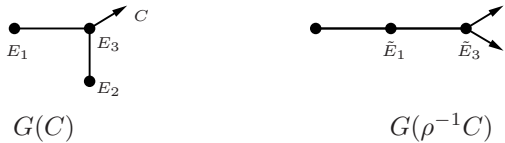
\begin{center}
\begin{texdraw}
\arrowheadtype t:F \arrowheadsize l:0.08 w:0.04 \drawdim mm \setgray
0 \move(10 20) \fcir f:0 r:0.8 \rlvec (10 0) \fcir f:0 r:0.8 \rlvec
(0 -7) \fcir f:0 r:0.8 \move (20 20) \avec (25 23) \move (26 21)
\htext{\tiny{$C$}}

\move (9 16) \htext{\tiny$E_1$} \move (21 17) \htext{\tiny{$E_3$}}
\move (21 11) \htext{\tiny{$E_2$}}

\move(50 20)
\fcir f:0 r:0.8 \rlvec (10 0) \fcir f:0 r:0.8 \rlvec (10 0) \fcir
f:0 r:0.8 \avec (75 23) \move (70 20) \avec (75 17)

\move(59 16) \htext{\tiny{$\tilde{E}_1$}} \move(69 16)
\htext{\tiny{$\tilde{E}_3$}}

\move (10 5) \htext{$G(C)$} \move (60 5) \htext{$G(\rho^{-1}C)$}
\end{texdraw}
\end{center}
where $\tilde{E}_1$, $\tilde{E}_3$ are the vertices associated to
$E_1$ and $E_3$ respectively.

Consider now a curve $C$ with characteristic exponents $\{4,6,7\}$.
Take  $\rho$ the ramification given by $\rho(u,v)=(u^4,v)$ and put
$\tilde{C}=\rho^{-1}C$. Then we have that
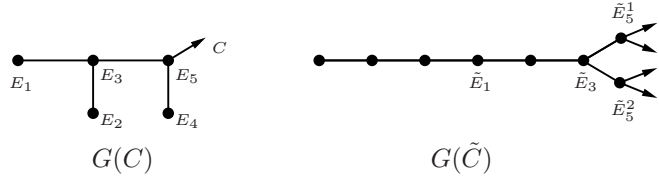
\vspace{\baselineskip}
\begin{center}
\begin{texdraw}
\arrowheadtype t:F \arrowheadsize l:0.08 w:0.04 \drawdim mm \setgray
0 \move(10 20) \fcir f:0 r:0.8 \rlvec (10 0) \fcir f:0 r:0.8 \rlvec
(0 -7) \fcir f:0 r:0.8 \move (20 20) \rlvec (10 0) \fcir f:0 r:0.8
\rlvec (0 -7) \fcir f:0 r:0.8 \move (30 20) \avec (35 23)

\move (36 21) \htext{\tiny{$C$}}

\move (9 16) \htext{\tiny$E_1$} \move (21 11) \htext{\tiny{$E_2$}}
\move (21 17) \htext{\tiny{$E_3$}} \move (31 11)
\htext{\tiny{$E_4$}} \move (31 17) \htext{\tiny{$E_5$}}

\move(50 20)
\fcir f:0 r:0.8 \rlvec (7 0) \fcir f:0 r:0.8 \rlvec (7 0) \fcir f:0
r:0.8 \rlvec (7 0) \fcir f:0 r:0.8 \rlvec (7 0) \fcir f:0 r:0.8
\rlvec (7 0) \fcir f:0 r:0.8 \rlvec (5 3) \fcir f:0 r:0.8 \avec (95
25)  \move (90 23) \avec (95 21)

\move (85 20) \rlvec (5 -3) \fcir f:0 r:0.8

\move (90 17) \avec (95 19) \move (90 17) \avec (95 15)

\move(70 16) \htext{\tiny{$\tilde{E}_1$}} \move(84 16)
\htext{\tiny{$\tilde{E}_3$}}

\move (89 25) \htext{\tiny{$\tilde{E}_{5}^{1}$}} \move (89 12)
\htext{\tiny{$\tilde{E}_{5}^{2}$}} \move (20 5) \htext{$G(C)$} \move
(65 5) \htext{$G(\tilde{C})$}
\end{texdraw}
\end{center}
Note that $E_3$ has one associated vertex $\tilde{E}_3$  and that
$E_5$ has two associated vertices $\tilde{E}_{5}^{1}$ and
$\tilde{E}_{5}^{2}$ in $G(\tilde{C})$.
\end{Example}

\begin{Remark}\label{nota:E1-bifurcacion}
Let us denote by $\tilde{E}_{1}$  the divisor of $G(\tilde{C})$ with
$v(\tilde{E}_{1}) =n$. It is unique since it  precedes all the other
bifurcation divisors and it could be or not a bifurcation divisor.
Moreover, $\tilde{E}_1$ is a bifurcation divisor of $G(\tilde{C})$
if and only if $E_1$ is a bifurcation divisor of $G(C)$ and
$b_{\tilde{E}_1}=b_{E_1}$. Then, the divisor $E_1$ of $G(C)$ has
always a unique divisor, denoted by $\tilde{E}_1$, which is
associated to it in $G(\tilde{C})$ even if $E_1 \not \in B(C)$.
Recall that $E_1$ is a bifurcation divisor if and only if the number
of different tangent lines in the tangent cone of $C$ is $\geq 2$.
\end{Remark}
We have seen that there is a bijection between the Puiseux series of
$C_i$ and the irreducible components of $\rho^{-1}C_i$. In
particular,  this implies that the choice of a vertex $\tilde{E}^l
\in B(\tilde{K}_i)$ associated to a bifurcation divisor $E$ is
equivalent to the choice of a $\underline{n}_{E}$-th root of the
unity $\xi_l$. Thus there are $e_{E}^{i}= n^{i}/\underline{n}_E$
irreducible components $\{\sigma_{lt}^{i}\}_{t=1}^{e_{E}^{i}}$ of
$\rho^{-1}C_i$ such that $\tilde{E}^l$ belongs to their geodesics.
Moreover, the curve $\sigma_{lt}^{i}$ is given by
$(v-\eta_{lt}^{i}(u)=0)$ where
$$
\eta_{lt}^{i}(u)=\sum_{s \geq n^i} a_{s}^{i} (\zeta_{ilt})^s
u^{sn/n^i}, \text{ for } t=1, \ldots, e_{E}^{i}.
$$
and $\{\zeta_{ilt}\}_{t=1}^{e_{E}^{i}}$ are the $e_{E}^{i}$-th roots
of $\xi_l$. Additionally , if $\gamma_E$ is an $E$-curvette of a
bifurcation divisor $E$ of $G(C)$, the curve $\rho^{-1}\gamma_E$ has
$m(E)=\underline{n}_E n_E$ irreducible components which are all
non-singular and each divisor $\tilde{E}^l$ belongs to the geodesic
of exactly $n_E$ branches of $\rho^{-1}\gamma_E$ which are curvettes
of $\tilde{E}^l$ in different points. In particular, we can prove
the following result

\begin{Lemma}\label{lema:apB}
Let $E$ be either a  bifurcation divisor of $G(C)$  or $E=E_1$ and
consider any of its associated divisors $\tilde{E}$ in
$G(\tilde{C})$. Then there exists a morphism $\rho_{\tilde{E},E} :
\tilde{E}_{red} \rightarrow E_{red}$ which is a ramification of
order $n_E$.
\end{Lemma}
\begin{proof}
Consider a $C$-ramification $\rho: ({\mathbb C}^2,0) \rightarrow
({\mathbb C}^2,0)$ given by $x=u^n, y=v$.  Let $\pi_{\tilde{E}}:
\tilde{M}_{\tilde{E}} \rightarrow ({\mathbb C}^2,0)$ be the
reduction of $\pi_{\tilde{C}}$ to $\tilde{E}$ and $\pi_{E}: M_E
\rightarrow ({\mathbb C}^2,0)$ be the reduction of $\pi_C$ to $E$.
Let us define the map $\rho_{\tilde{E},E} : \tilde{E}_{red}
\rightarrow E_{red}$. The map $\rho_{\tilde{E},E}$ sends the
``infinity point" of $\tilde{E}_{red}$ (that is, the origin of the
second chart of $\tilde{E}_{red}$) into  the ``infinity point" of
$E_{red}$. For any other  point $P$ of $\tilde{E}_{red}$, we
consider an $\tilde{E}$-curvette
$\gamma_{\tilde{E}}^P=(v-\psi_{\tilde{E}}^{P}(u)=0)$ with
$$
\psi_{\tilde{E}}^{P}(u) = \sum_{i=1}^{v(\tilde{E})-1}
a_i^{\tilde{E}} u^{i} + a_{v(\tilde{E})}^{P} u^{v(\tilde{E})},
$$
and such that $\pi_{\tilde{E}}^* \gamma_{\tilde{E}}^P \cap
\tilde{E}_{red} = \{P\}$. Let $\gamma_{E}^P$ be the curve given by
the Puiseux series
$$
y^P(x)= \sum_{i=1}^{v(\tilde{E})-1} a_i^{\tilde{E}} x^{i/m(E)} +
a_{v(\tilde{E})}^{P} x^{v(\tilde{E})/m(E)}.
$$
Thus $\gamma_E^P$ is an $E$-curvette and we define
$\rho_{\tilde{E},E}(P)$ to be the only point $\pi_{E}^* \gamma_{E}^P
\cap E_{red}$. From the properties of $\rho$ we deduce that
$\rho_{\tilde{E},E}$ is a ramification of order $n_E$.
\end{proof}
Remark also that, if $\gamma_{E_t}$ is a curvette of a terminal
divisor $E_t$ of a dead arc with bifurcation divisor $E$, then
$\rho^{-1}\gamma_{E_t}$ is composed by $m(E_t)=\underline{n}_{E}$
non-singular irreducible components and each divisor $\tilde{E}^l$
belongs to the geodesic of exactly  one branch of
$\rho^{-1}\gamma_{E_t}$, where
$\{\tilde{E}^l\}_{l=1}^{\underline{n}_E}$ are the divisors
associated to $E$ in $G(\tilde{C})$.

For more results concerning foliations, ramifications and blow-ups,
the reader can refer to \cite{For}.


\begin{thebibliography}{99}
\bibitem{Bri-K} E. Brieskorn; H. Kn\"{o}rrer: \textit{Plane algebraic curves.}
Birkh\"auser Verlag, Basel, 1986.

\bibitem{Cam-S} C. Camacho; P. Sad:
\textit{Invariant Varieties Through Singularities of Holomorphic
Vector Fields.} Ann. of  Math. (2) \textbf{115}, 3 (1982), 579--595.

\bibitem{Cam-S-LN}  C. Camacho; A. Lins Neto; P. Sad: \textit{Topological Invariants and Equidesingularisation for Holomorphic
Vector Fields.} J. Differential Geometry, \textbf{20}, 1 (1984),
143--174.

\bibitem{Cas-00} E. Casas-Alvero: \textit{Singularities of plane
curves.} London Mathematical Society Lecture Note Series, vol. 276,
Cambridge University Press, Cambridge, 2000.


\bibitem{Cor-03} N. Corral: \textit{Sur la topologie des courbes polaires de certains
feuilletages singuliers.} Ann. Inst. Fourier (Grenoble) \textbf{53}.
3 (2003), 787-814.

\bibitem{Cor-06} N. Corral: \textit{D\'etermination du type d'\'equisingularit\'e polaire.}
C. R. Math. Acad. Sci. Paris, S\'er. I, \textbf{344}, 1 (2007),
33--36.

\bibitem{Cor-08} N. Corral: \textit{Polar pencil of curves and
foliations.} To appear in Asterisque.

\bibitem{Fas-P} T. Fassarella; J. V. Pereira: \textit{On the degree
of polar transformations. An approach through logarithmic
foliations.} Sel. Math. (N. S.) \textbf{13}, 2 (2007), 239-252.

\bibitem{For} P. Fortuny Ayuso: \textit{Ramifications and Singularities of
Foliations.} {Proceedings of ``Transgressive Computing, a conference
in honour of Jean Della Dora"}, Granada (2006), 247--256.

\bibitem{Gar} E. R. Garc\'{\i}a Barroso: \textit{Sur les courbes polaires d'une courbe plane r\'eduite.}
Proc. London Math. Soc. (3) \textbf{81}, 1 (2000), 1--28.

\bibitem{Kuo-L} T. C. Kuo; Y. C. Lu: \textit{On analytic function germs of two complex
variables.} Topology \textbf{16}, 4 (1977), 299--310.

\bibitem{Le-M-W} {L\^{e} D\~{u}ng Tr\`{a}ng; F. Michel; C. Weber:}
\textit{Sur le comportement des polaires associ\'ees aux germes de
              courbes planes.} Compositio Math. \textbf{72}, 1
              (1989), 87--113.

\bibitem{Mer} M. Merle: \textit{Invariants polaires des courbes
planes.} Invent. Math. \textbf{41}, 2 (1977), 103--111.

\bibitem{Mol} R. Mol:  \textit{Classes polaires associ\'{e}s aux
distributions holomorphes de sous-espaces tangents.} Bull. Braz.
Math. Soc. (N. S.), \textbf{37}, 1 (2006), 29-48.

\bibitem{Rou} P. Rouill\'{e}: \textit{Th\'eor\`eme de {M}erle: cas des 1-formes de type courbes
              g\'en\'eralis\'ees.} Bol. Soc. Brasil. Mat. (N.S.)
              \textbf{30}, 3 (1999), 293--314.

\bibitem{Tei} B. Teissier:  \textit{Vari\'et\'es polaires I.} Invent. Math. {\bf 40}, (1977),
267-292.

\bibitem{Zar} O. Zariski: \textit{General theory of saturation and
               of saturated local rings II.} Amer. J. Math.
               \textbf{93} (1971), 872--964.

\bibitem{Zar2} O. Zariski: \textit{Studies in equisingularity I, II, III.}
Amer. J. of Math. \textbf{87}, 2 (1965), 507--533; \textbf{87}, 4
(1965), 972--1006; \textbf{90}, 3 (1968), 961--1023.

\end{thebibliography}
\end{document}

